\documentclass[letterpaper, 10 pt, conference]{ieeeconf}  

\IEEEoverridecommandlockouts  

\usepackage{cite}
\usepackage{amsmath,amssymb,amsfonts}
\usepackage{textcomp}

\def\BibTeX{{\rm B\kern-.05em{\sc i\kern-.025em b}\kern-.08em
    T\kern-.1667em\lower.7ex\hbox{E}\kern-.125emX}}
\usepackage{graphicx}
\usepackage{bm}
\usepackage{mathtools}
\usepackage{amsmath,amsfonts,amssymb}
\usepackage{algorithm,algpseudocode}
\usepackage{url}
\usepackage{cases}
\usepackage{ushort}

\DeclareMathOperator{\col}{col}

\begin{document}

\title{Two-Facet Scalable Cooperative Optimization of Multi-Agent Systems in The Networked Environment}
\author{X. Huo, \IEEEmembership{Student Member, IEEE}, M. Liu \IEEEmembership{Member, IEEE}
\thanks{X. Huo and M. Liu are with the Department of Electrical and Computer Engineering, University of Utah, Salt Lake City, UT 84112 USA (e-mail: xiang.huo, mingxi.liu@utah.edu).}}

\maketitle

\begin{abstract}
Cooperatively optimizing a vast number of agents that are connected over a large-scale network brings unprecedented scalability challenges. This paper revolves around problems optimizing coupled objective functions under coupled network-induced constraints and local constraints. The scalability of existing optimization paradigms is limited by either the agent population size or the network dimension. As a radical improvement, this paper for the first time constructs a two-facet scalable decentralized optimization framework. To this end, we first develop a systemic network dimension reduction technique to virtually cluster the agents and lower the dimension of network-induced constraints, then constitute a novel shrunken-primal-multi-dual subgradient (SPMDS) algorithm based on the reduced-dimension network. Rigorous optimality and convergence analyses of the proposed decentralized optimization framework are provided. The SPMDS-based optimization framework is free of agent-to-agent communication and no additional aggregators are required for agent clusters. The efficiency and efficacy of the proposed approaches are demonstrated, in comparison with benchmark methods, through simulations of electric vehicle charging control problems and traffic congestion control problems.

\end{abstract}

\section{Background}
\label{sec:introduction}
Cooperatively controlling and optimizing multiple agents in a networked computing environment has gained growing attention in the recent research \cite{rahbari2014incremental, liu2017decentralized,morrisse2012distributed,yang2013consensus,he2007towards}. These problems, referred to as networked optimization problems (NOP), aim to coordinate the agents to achieve a system-level goal while considering agents' collective impacts on the network environment. One popular NOP is the resource allocation problem which deals with allocating available resources to a number of agents. A generic NOP formulation can be written as \cite{xu2017distributed}
\begin{subequations} \label{eq:1snn}
\begin{align}
    &\underset{\bm{x}}{\text{min}} \quad {\sum_{i=1}^{v}f_i(\bm{x}_i)} \label{eq:1snna}\\
&\text{s.t.}  \ \quad \bm{x}_{i} \in \mathbb{X}_{i}, \quad \forall i=1,2, \ldots, v, \label{eq:1snnb}\\
& \,\,\, \qquad \bm{C}\bm{x} = \bm{d}, \label{eq:1snnc}
\end{align}
\end{subequations}
where $\bm{x}=\col(\bm{x}_1,\ldots,\bm{x}_v)$ with $\col(\cdot)$ denoting a column vector and each $\bm{x}_i$ denoting the resources allocated to agent $i$ associated with a local cost function $f_i(\bm{x}_i)$, and $\mathbb{X}_{i}$ denotes the local constraint set for agent $i$. Eqn. \eqref{eq:1snnc} represents the networked resource allocation constraint where $\bm{C}$ denotes a non-negative matrix with only one non-zero entry in each column and $\bm{d}$ is a given vector. Then, the resource allocation problem in \eqref{eq:1snn} aims at finding the minimum total cost under local constraints and global equality constraints. Problem in \eqref{eq:1snn} has a wide range of applications, including multiple resources allocation in energy management systems \cite{rahbari2014incremental}, power regulation \cite{morrisse2012distributed}, economic dispatch \cite{yang2013consensus}, etc.

Another NOP that has received significant attention is the cooperative control of multiple agents. One popular application is the electric vehicle (EV) charging control problem \cite{liu2017decentralized}. The charging processes of a large number of EVs require coordination and control for the alleviation of impacts on the distribution network (networked constraint) and for the provision of various grid services (system-level objective) \cite{lopes2010integration,clement2009impact,bronzini2011coordination,qian2010modeling}. A generic EV charging control problem can be written as \cite{liu2017decentralized}
\begin{subequations} \label{eq:1sn}
\begin{align}
    &\underset{\bm{x}}{\text{min}} \quad {F(\bm{x})}+{\sum_{i=1}^{v}f_i(\bm{x}_i)} \label{eq:1sa}\\
&\text{s.t.}  \ \quad \bm{x}_{i} \in \mathbb{X}_{i}, \quad \forall i=1,2, \ldots, v, \label{eq:1sb}\\
& \,\,\, \qquad \bm{A}\bm{x} \leq \bm{b}, \label{eq:1sc}
\end{align}
\end{subequations}
where $F(\bm{x})$ denotes the grid service objective, e.g., valley-filling and power trajectory tracking, $f_i(\bm{x}_i)$ represents local interests, $\mathbb{X}_i$ denotes local charging physical limits, and \eqref{eq:1sc} denotes distribution network constraints, e.g., nodal voltage magnitudes must be maintained within a certain range. Unlike the problem in \eqref{eq:1snn}, $\bm{x}_i$ in the objective function $F(\bm{x})$ in \eqref{eq:1sa} is normally strongly coupled and non-separable, and a more general global inequality constraint \eqref{eq:1sc} is in place.

In contrast to the problem in \eqref{eq:1snn}, the problem in \eqref{eq:1sn} can cover a broader range of applications. Besides EV charging control, it can be used for rate control in communication networks \cite{kelly1998rate} where the objective function is separable and denoted as the summation of the agents' utility functions \eqref{eq:1snna}, and subject to a global inequality constraint which represents the networked resource capacity limits. It also can be used for traffic congestion optimization in intelligent transportation systems where the objective function is non-separable due to the congestion cost \cite{he2007towards}.  In this paper, we are set to study efficient algorithms for the generic formulation in \eqref{eq:1sn} that has a wide range of network optimization applications. To make the discussions more practical, we will elicit the problem formulation from EV charging control problems.


\section{Related Work}
A variety of optimization paradigms have been researched for \eqref{eq:1sn}, including centralized, distributed, decentralized, and hybrid optimization \cite{guo1996algorithm,richardson2011optimal,yang2016distributed,siljak2011decentralized,wang2017hybrid}. Centralized strategies \cite{kang2015centralized,jian2017high,lozano2004centralized}, though being theoretically feasible and efficient, can hardly be deployed in large-scale applications mainly due to the impaired scalability \emph{w.r.t.} both the agent population size and network dimension. From the practicability perspective, centralized algorithms require private information sharing from agents to a central operator (CO) which poses potential risks to customer privacy \cite{yang2014improved}; they also lack robustness, e.g., a small change in the network may lead to the redesign of the whole centralized algorithm \cite{yang2013consensus}. 

The deficiency of centralized optimization has largely driven the development of distributed optimization.  In \cite{de2019distrac}, a distributed low-communication protocol was designed for traffic congestion control based on vehicle-to-vehicle communication. Rahman \emph{et al.}\cite{rahman2019vehicle} developed a switchable optimization-incorporated distributed controller that utilizes a sparse communication network. More recently, the alternating direction method of multipliers (ADMM) \cite{boyd2011distributed} has become popular in distributed optimization problems that aim at the minimization of separable cost-related objective functions \cite{du2019admm,peng2014distributed,rivera2016distributed}. The abovementioned approaches, together with other existing distributed methods, indeed alleviate the scalability issue, however, the inevitable peer-to-peer communication poses cyber-security risks to the entire control framework. To eliminate the frequent peer-to-peer communication and retain the scalability, decentralized (or hybrid) optimization for separable objective functions has attracted more attentions \cite{braun2016distributed,wang2017hybrid,wen2012decentralized}. A hybrid centralized-decentralized EV charging control framework developed in \cite{wang2017hybrid} is able to eliminate the interactions between EVs. A hierarchical, iterative distributed optimization algorithm was developed in  \cite{braun2016distributed}, which does not require communication between agents. Similarly, \cite{wen2012decentralized} mitigates the security problems inherently in EV charging control problem by utilizing a decentralized charging selection concept.

However, assuming separability of decision variables in the objective function and constraints appears impractical in more general cases. For example, in EV charging control for valley-filling or power trajectory tracking, all EVs' decision variables are strongly coupled and non-separable in the objective function and constraints. In this case, ADMM-based algorithms cannot work effectively, so do many other distributed or decentralized algorithms, e.g. \cite{xu2017distributed,yang2013consensus,johansson2009distributed}. In the state of the art, few theoretical results have attempted to address this issue. In \cite{koshal2011multiuser}, a regularized primal dual subgradient (RPDS) algorithm was developed by approximating the difference between the original optimal solution  and the regularized counterpart. However, the regularization that was introduced to guarantee convergence led to inevitable convergence errors. As an improvement, Liu \emph{et al.} in \cite{liu2017decentralized} developed a shrunken primal dual subgradient (SPDS) algorithm whose convergence does not rely on the regularization and eliminates convergence errors. Admittedly, RPDS and SPDS are more general and are both scalable  \emph{w.r.t.} the number of agents, however, the network dimension presents a hidden scalability issue. Specifically, computations in the algorithm iterations must involve either the network connectivity matrix or the adjacent matrix whose dimensions dramatically increase as the distribution network dimension grows. Hence, directly implementing RPDS and SPDS in large-scale distribution networks without considering the network dimension will cause memory overflow and may exceed the computing capacity of the on-board controller. 

In order to overcome the scalability issue induced by the network dimension, network division or approximation strategies have been widely investigated. In \cite{zhou2020gradient}, a large distribution system was divided into small areas to improve the convergence speed for distribution system state estimation. Similarly, \cite{du2019admm} partitioned the network into regional subsystems, then developed an ADMM-based distributed state estimation algorithm to resolve cyber attacks within the subsystems. However, peer-to-peer communication prevalently exists in both \cite{zhou2020gradient} and \cite{du2019admm} which may cause  potentially cyber-security risks.  Alternatively, an approximation methodology was proposed in \cite{moruarescu2012dimension} to transform a large-scale networked system into a lower dimensional one, however, the inexact approximation introduces inevitable precision deterioration. In our preliminary work \cite{2020cdcspmds}, a shrunken primal-multi-dual subgradient (SPMDS) was developed based on a heuristic dimension reduction strategy. Unfortunately, the developed strategy is not systemic and theoretical convergence and optimality guarantees of SPMDS are missing. More importantly, existing dimension reduction methods are incapable of being integrated into existing decentralized optimization frameworks without sacrificing convergence and optimality. Considering the scalability gaps in the NOP algorithm research, this paper is set to develop a decentralized algorithm that has two-facet scalability and provide convergence and optimality guarantees. 

The contribution of this paper is four-fold: (1) This paper, for the first time, designs a novel systemic dimension reduction strategy that partitions the network and groups primal decision variables to reduce the network complexity; (2) A novel decentralized optimization framework that synthesizes SPMDS and the dimension reduction technique is developed to achieve the scalability  \emph{w.r.t.}  both the number of primal decision variables and the network dimension; (3) The proposed framework does not impose additional communication load, i.e., no additional aggregators are required for the virtual agent groups, and all dual variables can be updated in parallel; (4) Optimality conditions, convergence guarantees, and complexity relaxation of SPMDS in conjunction with the dimension reduction technique are, for the first time, rigorously analyzed and proved.

The remainder of this paper is organised as follows: In Section III, we elicit the strongly coupled NOP from the EV charging control problem. Section IV presents the main results of this paper, including the dimension reduction strategy and SPMDS. Convergence analysis and computational complexity analysis of the proposed approaches are also discussed in Section IV. We give simulation results of different scenarios, in comparison with benchmark methods, in Section V to show the efficacy and efficiency of the proposed approaches. Section VI concludes this paper.

\section{Preliminaries and Problem Formulation}

In this section, we build a radial power distribution network model, based on which we formulate an EV charging control problem which is a representative large-scale NOP with non-separable objective function and coupled constraints that well manifests the problem in \eqref{eq:1sn}. The purpose of introducing the EV charging control problem is to offer a solid practical background of the approaches proposed in Section IV, while the applicability of the approaches is not necessarily confined with this particular application as shown in Section V.

\subsection{Radial Distribution Network Model}

Power flow of a radial distribution network can be described by the DistFlow branch equations which only involve the real power load, reactive power load, and voltage magnitude \cite{baran1989network}. Let $\mathbb{N} = \{0,1,\ldots,n \}$ denote the set of nodes in a radial distribution network, where Node 0 is the slack node that maintains its voltage magnitude at a constant $V_0$, $\mathbb{C}_j$ denotes the set of bus $j$'s children. For two adjacent nodes namely Node $i$ and Node $j$, let $(i,j)$ denote the line segment connecting them, $\mathcal{P}_{ij}$ and $\mathcal{Q}_{ij}$ denote the active and reactive power flow from Node $i$ to Node $j$ respectively, and $r_{ij}$ and $x_{ij}$ denote the resistance and reactance of line $(i,j)$, respectively. For each Node $j$, let $P_j$ and $Q_j$ denote its active power and reactive power consumption, respectively, and let $p_j$ and $q_j$ denote its active and reactive power injection, respectively. 
\begin{figure}[!htb]
    \centering
    \includegraphics[width=0.45\textwidth]{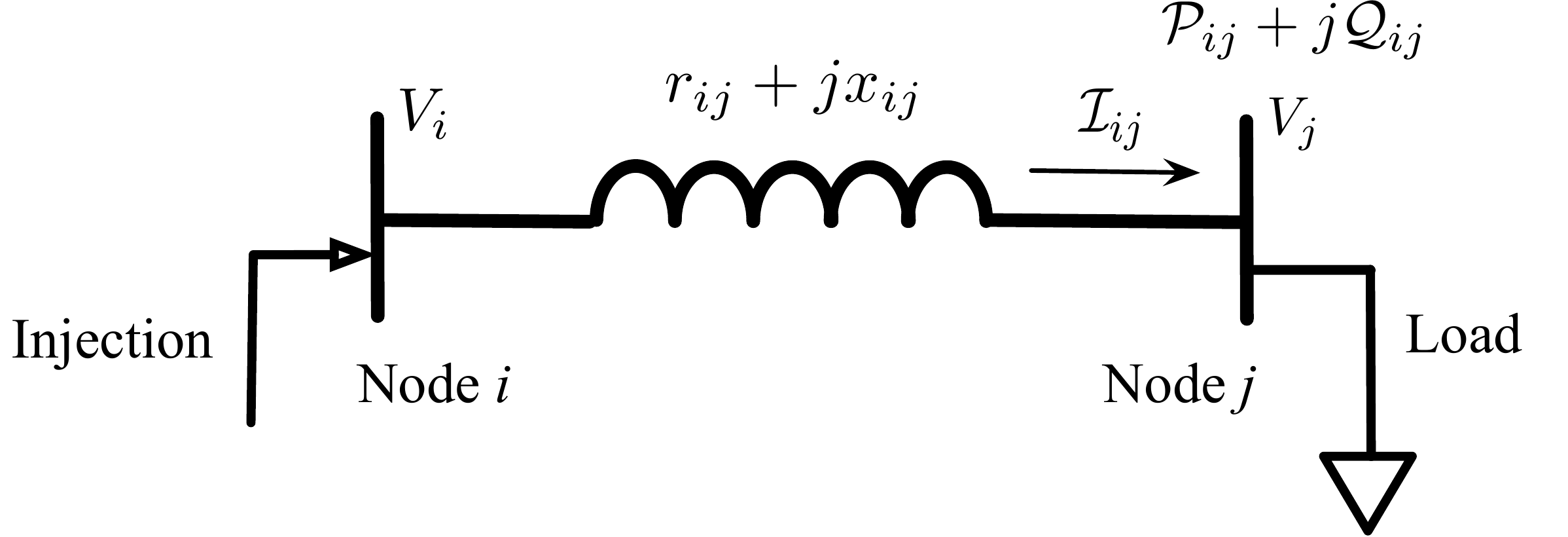}
    \caption{A basic two node power flow model \cite{shigenobu2017demand}}
    \label{twobus}
\end{figure} 
A basic two node power flow model is shown in Fig. \ref{twobus}, and the power flow of the radial distribution network can be defined through the DistFlow branch equations \cite{baran1989network} as 
\begin{subequations}
\begin{align}
    \mathcal{P}_{ij} - \sum_{k \in \mathbb{C}_j} \mathcal{P}_{jk} &= P_j  - p_j + r_{ij}\mathcal{I}_{ij}^2, \\
    \mathcal{Q}_{ij} - \sum_{k \in \mathbb{C}_j} \mathcal{Q}_{jk} &= Q_j  - q_j + x_{ij}\mathcal{I}_{ij}^2, \\
    V_i^2 - V_j^2 &= 2(r_{ij} \mathcal{P}_{ij} + x_{ij}\mathcal{Q}_{ij}){-}(r_{ij}^2 + x_{ij}^2)\mathcal{I}_{ij}^2,
\end{align}
\end{subequations}
where $\mathcal{I}_{ij}^2 = (\mathcal{P}_{ij}^2+\mathcal{Q}_{ij}^2)/V_i^2$.

To simplify the network model, a DistFlow model can be linearized to the LinDistFlow model by omitting the line loss and some higher order terms, i.e., $r_{ij}\mathcal{I}_{ij}^2, x_{ij}\mathcal{I}_{ij}^2 \ \text{and} \ (r_{ij}^2 + x_{ij}^2)\mathcal{I}_{ij}^2$ \cite{baran1989optimal}. It has been shown in \cite{bansal2014plug,gan2014convex,zhang2016scalable} that this linearization has negligible impacts on the model accuracy, by which the error introduced is relatively small and normally on the order of $1\%$ \cite{hu2020voltage,zhu2015fast}. Hence, in this paper, we adopt the LinDistFlow model to simplify the power flow description and better illustrate the algorithm design. The LinDistFlow model which maps nodal real and reactive power consumption can be represented as \cite{baran1989optimal}
\begin{equation}
\bm{V}(T)=\bm{V}_{0}-2 \bm{R} \bm{P}(T)-2 \bm{X} \bm{Q}(T),
\label{eq:1}
\end{equation}
where $\bm{V}(T) \in \mathbb{R}^{n}$ consist of the squared voltage magnitudes of Nodes 1 to $n$, $\bm{V}_0 = V_0^2 \bm{1}_n \in \mathbb{R}^{n}$ denotes the slack constant voltage magnitude vector, and $\bm{P}(T) \in \mathbb{R}^{n}$ and $\bm{Q}(T) \in \mathbb{R}^{n}$ denote the real and reactive power consumption from Node 1 to Node $n$, respectively. In addition, $\bm{R}$ and $\bm{X}$ are the voltage-to-power-consumption sensitivity matrices \cite{zhou2020gradient} defined as
\begin{equation} \label{3}
\begin{aligned}
\bm{R} &\in \mathbb{R}^{n \times n}, \quad R_{ij}=\sum_{(i, j) \in \mathbb{E}_{i} \cap \mathbb{E}_{j}} 
r_{ij}, \\
\bm{X}  &\in \mathbb{R}^{n \times n}, \quad 
X_{ij}=\sum_{(i,j) \in \mathbb{E}_{i} \cap \mathbb{E}_{j}} 
x_{ij},
\end{aligned}
\end{equation}
where $\mathbb{E}_{i}$ and $\mathbb{E}_{j}$ are the line sets connecting Node 0 and Node $i$, and Node 0 and Node $j$, respectively \cite{farivar2013equilibrium}, i.e., the voltage-to-power-consumption sensitivity factors $R_{ij} (X_{ij})$ are obtained through the resistance (reactance) of the common path of  Node $i$ and Node $j$ leading back to Node 0 \cite{zhou2020gradient}.

\subsection{EV Charging in The Distribution Network}
In this paper, we consider that the load at each node is composed of uncontrollable baseline load (e.g., lighting, A/C) and EV charging load. Without the loss of generality, we assume that the EV charging load is the only controllable load and EVs only consume real power \cite{liu2017decentralized}. We also assume that the baseline load can be perfectly forecasted. Uncertainties will be investigated in our future work. Consequently, the LinDistFlow model in \eqref{eq:1} can be rewritten as 
\begin{equation}
\bm{V}(T)=\bm{V}_{0}-\bm{V}_{b}(T)-2 \bm{R} \bm{G}\Bar{\bm{P}}\bm{u}(T),
\label{eq:5}
\end{equation}
where $\bm{V}_{b}(T)$ denotes the voltage drop caused by the baseline load at time $T$, $\bm{u}(T){=}\col(u_1(T)\cdots u_v(T)) \in \mathbb{R}^{v}$ in $[\bm{0}, \bm{1}]$ contains the normalized charging rates of all $v$ EVs connected at the distribution network, $\bm{G}$ is the charging aggregation matrix that aggregates the charging power of EVs connected at the same node, and $\bm{\Bar{P}}$ is the maximum charging power matrix of all EVs. Herein, $\bm{G}$ and $\bm{\Bar{P}}$ are defined by 
\begin{equation}
\bm{G} \triangleq \bigoplus_{j = 1}^{n} \bm{G}_{j} \in \mathbb{R}^{n \times v},~\bm{\Bar{P}} \triangleq \bigoplus_{i = 1}^{v} \Bar{P}_{i} \in \mathbb{R}^{v \times v}, \nonumber
\end{equation}
where $\bigoplus$ denotes the matrix direct sum hereinafter, $\bm{G}_j = \bm{1}_{v_j}^\mathsf{T}$ is the charging power aggregation vector, ${\Bar{P}}_i$ is the maximum charging power of the $i$th EV, and $v_j$ is the number of EVs connected at node $j$ with $\sum_{j=1}^{n}v_j=v$. To simplify the presentation of \eqref{eq:5}, let $\bm{V}_c(T)$ denote $\bm{V}_{0}-\bm{V}_{b}(T)$ and $\bm{D} \in \mathbb{R}^{n \times v}$ denote $2 \bm{R} \bm{G}\Bar{\bm{P}}$, then we have
\begin{equation}
\bm{V}(T)=\bm{V}_c(T) - \bm{D}\bm{u}(T).
\label{eq:6}
\end{equation}

Assume the valley-filling service duration is fixed with $K$ time intervals, where $k$ and $k+K-1$ denote the starting and ending time, respectively. By augmenting the system output in \eqref{eq:6} along the valley-filling period $[k,k+K-1]$, we have
\begin{equation}
   \bm{V}(k) = \bm{V}_{c}(k) - \sum_{i=1}^{v} \mathcal{D}_{i} \mathcal{U}_{i}(k),
   \label{eq:8}
\end{equation}
where 
\begin{align*}
 \bm{V}(k) &= \col\left(\bm{V}(k | k),\bm{V}(k+1 | k),
    \cdots, 
    \bm{V}(k+K-1 | k)
    \right), \\
\bm{V}_c{(k)}&=\col\left(
    \bm{V}_c(k), \bm{V}_c(k+1), \cdots, \bm{V}_c(k+K-1)
    \right), \\
 \mathcal{U}_i{(k)} &= \col\left(
     u_i(k | k),
     u_i(k+1 | k), 
    \cdots, 
    u_i(k+K-1|k)
    \right),\\
\mathcal{D}_{i} &= D_{i} \oplus \cdots \oplus D_{i} {\in} \mathbb{R}^{n K \times K},
\end{align*}
and $D_{i}$ denotes the $i$th column of $\bm{D}$. Time index $(k+j | k), \ j=0,\ldots,K-1$ denote the input (output) prediction at time $k+j$ based on the past knowledge up to time $k$ \cite{camacho2013model}.


Let $x_{i}(T)$ denote the energy remained to be charged to the $i$th EV at time $T$ and $B_{i}=-\eta_{i} \Delta t \bar{P}_{i}$ denote the maximum charging energy during time $\Delta t$ where $\eta_{i}$ is the charging efficiency, then the charging dynamics of the $i$th EV can be written as 
\begin{equation}
    x_{i}(T+1)=x_{i}(T)+B_{i} u_{i}(T).
\end{equation}
To guarantee all EVs are fully charged by the end of valley-filling, the following equality constraint must be satisfied
\begin{equation}
   \bm{x}(k)+\sum_{i=1}^{v} \mathcal{B}_{i,l} \mathcal{U}_{i}(k) = \bm{0},
   \label{eq:10}
\end{equation}
where  $\bm{x}(k) = \col\left(x_1(k)~ x_2(k)~\cdots~ x_v(k) \right) \in \mathbb{R}^{v}$, $\mathcal{B}_{i, l}=\left[B_{i, c}~ B_{i, c}~ \cdots~ B_{i, c}\right] \in \mathbb{R}^{n \times K}$, and $B_{i, c}$ denotes the $i$th column of the matrix $\bm{B} = \bigoplus_{i = 1}^{v}B_i $.

\subsection{Problem Formulation}
Valley-filling relies on the aggregated EV charging load to flatten the total demand curve of the distribution network. This can be achieved by minimizing the $\ell_2$-norm of the aggregated demand profile \cite{gan2012optimal}. By adopting the same notations as in \cite{liu2017decentralized} and dropping the time indicator $k$ in $\bm{V}(k)$, $\bm{V}_c(k)$, and $\mathcal{U}_i(k)$ hereinafter, we write the objective function of the valley-filling problem as 
\begin{align} 
\mathcal{F}(\mathcal{U}) = \frac{1}{2}\left\|P_{b}+\tilde{P} \mathcal{U}\right\|_{2}^{2} + \frac{\rho}{2} \| \mathcal{U} \|_{2}^{2},
\label{eq:11}
\end{align}
where $\mathcal{U} =\col(\mathcal{U}_1, \ldots, \mathcal{U}_v)\in \mathbb{R}^{vK}$, $P_{b} \in \mathbb{R}^K$ is the aggregated  baseline load profile along the valley-filling period, and $\tilde{P} \in \mathbb{R}^{K \times vK}$ is the aggregation matrix for all EVs' charging profiles. The last term in \eqref{eq:11} is a proxy of battery degradation cost which is approximated by a quadratic term of the charging rates \cite{ma2015distributed}.

The constraints of the valley-filling problem can be categorized into local constraints and networked constraints. For the local constraints, each EV should satisfy
\begin{equation}
  \mathcal{U}_{i} \in \mathbb{U}_{i},
    \mathbb{U}_{i} \triangleq \left\{\mathcal{U}_{i} | \mathbf{0} \leq \mathcal{U}_{i} \leq \mathbf{1}, x_{i}(k)+\mathcal{B}_{i,l} \mathcal{U}_{i}=0\right\}.
    \label{eq:12b}
\end{equation}
The purpose of \eqref{eq:12b} is to guarantee that each EV would be fully charged at the end of the valley-filling. For the distribution network, voltage magnitudes at all nodes should be limited within $\left[\underline{v}V_{0}, \bar{v}V_{0}\right]$, where $\underline{v}$ and $\overline{v}$ denote the lower and upper bounds, respectively. Let $\bm{V}_0=V_0^2\bm{1}_{nK}$, we have
\begin{equation}
\underline{v}^{2} \bm{V}_{0} \leq \bm{V} \leq \bar{v}^{2} \bm{V}_{0}.\label{eq:13}
\end{equation}
Then by using \eqref{eq:8} and with EVs' charging as the only controllable load, we only need to consider the lower bound for the network voltage constraint. This implies 
\begin{equation}
\bm{V}_c - \sum_{i=1}^{v} \mathcal{D}_{i} \mathcal{U}_{i} \geq \underline{v}^{2} \bm{V}_{0}. \label{eq:14}
\end{equation}
Note that, in more general cases where vehicle-to-grid is considered, the upper bound can be added back without affecting the algorithm design. Based on the above discussions, the valley-filling problem is formulated as
\begin{subequations} \label{eq:15}
\begin{align}
    &\underset{\mathcal{U}}{\text{min}} \quad {\mathcal{F}(\mathcal{U})} \label{eq:12sa}\\
&\text{s.t.}  \ \quad \mathcal{U}_{i} \in \mathbb{U}_{i}, \quad \forall i=1,2, \ldots, v, \label{eq:12sb}\\
& \,\,\, \qquad \underline{v}^{2} \bm{V}_{0} - \bm{V}_c + \sum_{i=1}^{v} \mathcal{D}_{i}\mathcal{U}_{i} \leq \bm{0}. 
\end{align}
\end{subequations}

\section{Main Results}
\subsection{Virtual Agent Grouping and Network Partitioning} \label{section_dimension_reduction}

For general linearized optimal power flow or EV charging control problems, decentralized optimization algorithms normally rely on the iterative primal and dual updates \cite{koshal2011multiuser,liu2017decentralized,yin2009nash}. State-of-the-art methods require the complete network topology information (e.g., connectivity matrix, adjacent matrix, or the sensitivity matrix $\bm{R}$ as in \eqref{3}) to execute the updates. For example, projection based algorithms including SPDS \cite{liu2017decentralized} and RPDS \cite{koshal2011multiuser} require complete network topology information, i.e., $\bm{A}$ matrix in \eqref{eq:1sc}. In the primal and dual updating processes, some interim matrices normally have the dimensions that are thousands of times of the full network dimension, leading to extra requirements on the on-board memory size and computing power. This severely impairs the scalability of those algorithms which were designed to be scalable. To overcome this, it is critical to reduce the dimension of the distribution network. 

The sensitivity matrix $\bm{R}$ maps the EV charging power to the distribution voltage profile $\bm{V}(T)$, e.g., each element $R_{ij}$ reflects the impact of the aggregated charging power at Node $j$ on the voltage magnitude of Node $i$. Since all $n$ nodes can have EVs connected, it is rational to consider the full dimension of $\bm{R}$ along the horizontal direction (all columns). However, the vertical direction (rows) of $\bm{R}$ can be reduced as EVs at a specific node (or some nodes as a cluster) only have \textit{major impacts} on a sub-vector of the original global voltage vector $\bm{V}(T)$.
In our preliminary work \cite{2020cdcspmds}, we proposed a heuristic way for dimension reduction, i.e., by concentrating and abstracting the major elements in the sensitivity matrix $\bm{R}$. Consequently, the elements in $\bm{R}$ which have relatively greater impacts on the nodal voltage magnitudes are preserved. However, this concentrating and abstracting process is complex, case-dependent, and not ready for the integration into decentralized algorithm design. In this paper, we propose a systemic dimension reduction approach that can be adopted for general network environment and seamlessly integrated into efficient decentralized algorithms.

\subsubsection{Virtually grouping EVs via $\mathcal{K}$-means}   

The $\mathcal{K}$-means clustering algorithm \cite{macqueen1967some} clusters data in a data set by minimizing the sum of squared distance. In analogy, EVs on nodes (agents) that have similar impacts on $\bm{V}(T)$ can be grouped. Based on $\mathcal{K}$-means, we partition the network nodes by solving
\begin{equation}
\underset{\bm{C}^{1},\ldots,\bm{C}^{r}}{\text{min}} \sum_{i=1}^{n}\underset{h=1,\ldots,r}{\text{min}}\left( \frac{1}{\kappa} \left\lVert [\bm{R}]_i - \bm{C}^h \right\rVert ^2_2 \right),
\label{3s}
\end{equation}
where $[\bm{R}]_i \in \mathbb{R}^n$ denotes the agent (node) $i$ which is the $i$th column of the sensitivity matrix $\bm{R}$, $\bm{C}^h$ is a column vector which denotes the nearest clustering center for $[\bm{R}]_i$, $r$ denotes the total number of desired clustering centers, and $\kappa$ denotes the number of agents (nodes). Given the number of clusters, the $\mathcal{K}$-means algorithm in \eqref{3s} finds the cluster centers $\bm{C}^{1},\ldots,\bm{C}^{r}$ such that
the sum-of-squared Euclidean distances from each point $[\bm{R}]_i$ to its nearest cluster center $\bm{C}^h$ are minimized.

Note that for the EV charging case, each column of the sensitivity matrix $\bm{R}$ represents a specific node. By applying the above proposed clustering strategy, all nodes are clustered based on the clustering results. Consequently, EVs connected at the corresponding nodes can be grouped accordingly. Therefore, this clustering strategy is independent of the number of EVs connected at each node, it is solely dependent on the physics of the network. An example of node clustering for a modified IEEE 13-bus test feeder used in \cite{liu2017decentralized} is shown in Fig. 1.
\begin{figure}[!htb]
    \centering
    \includegraphics[width=0.48\textwidth]{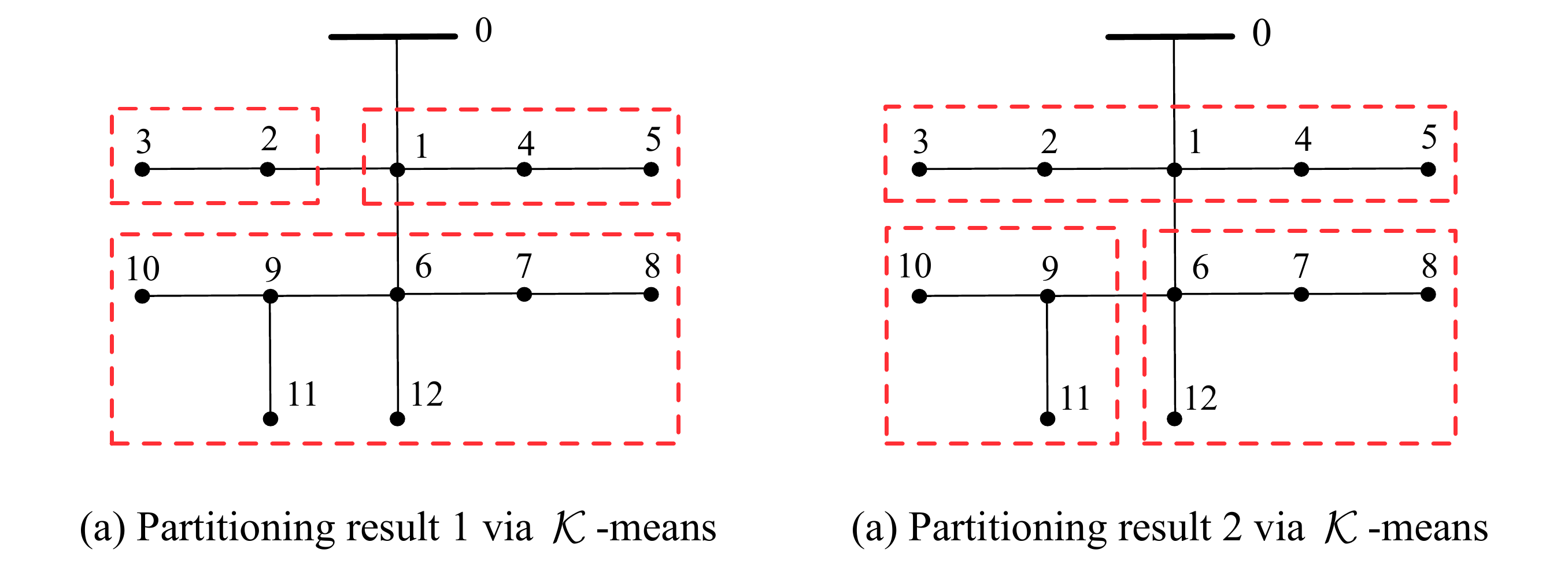}
    \caption{Network partitioning results via $\mathcal{K}$-means of the modified IEEE 13-bus test feeder (with clustering number $r=3$).}
    \label{fig:1}
\end{figure}

\noindent {\bf{Remark 1:}} $\mathcal{K}$-means
clusters $n$ data points into $r$ clusters so that the intra-group distances are low and inter-group distances are high, and $\mathcal{K}$-means converges until the clustering centers stop changing. The bilevel optimization programming problem in 
\eqref{3s} is NP-hard due to the 
nature of Euclidean sum-of-squares \cite{aloise2009np}. Therefore, the solution of $\mathcal{K}$-means converges to a local optimum without the guarantee of global optimality \cite{malinen2014balanced}. As discussed later in this section, the local optimum of clustering obtained through $\mathcal{K}$-means does not have any impact on the algorithm convergence. \hfill $\blacksquare$

\subsubsection{Determination of the voltage sub-vectors}  
\begin{figure}[!htb]
    \centering
    \includegraphics[width=0.4\textwidth, trim = 0mm 5mm 0mm 10mm, clip]{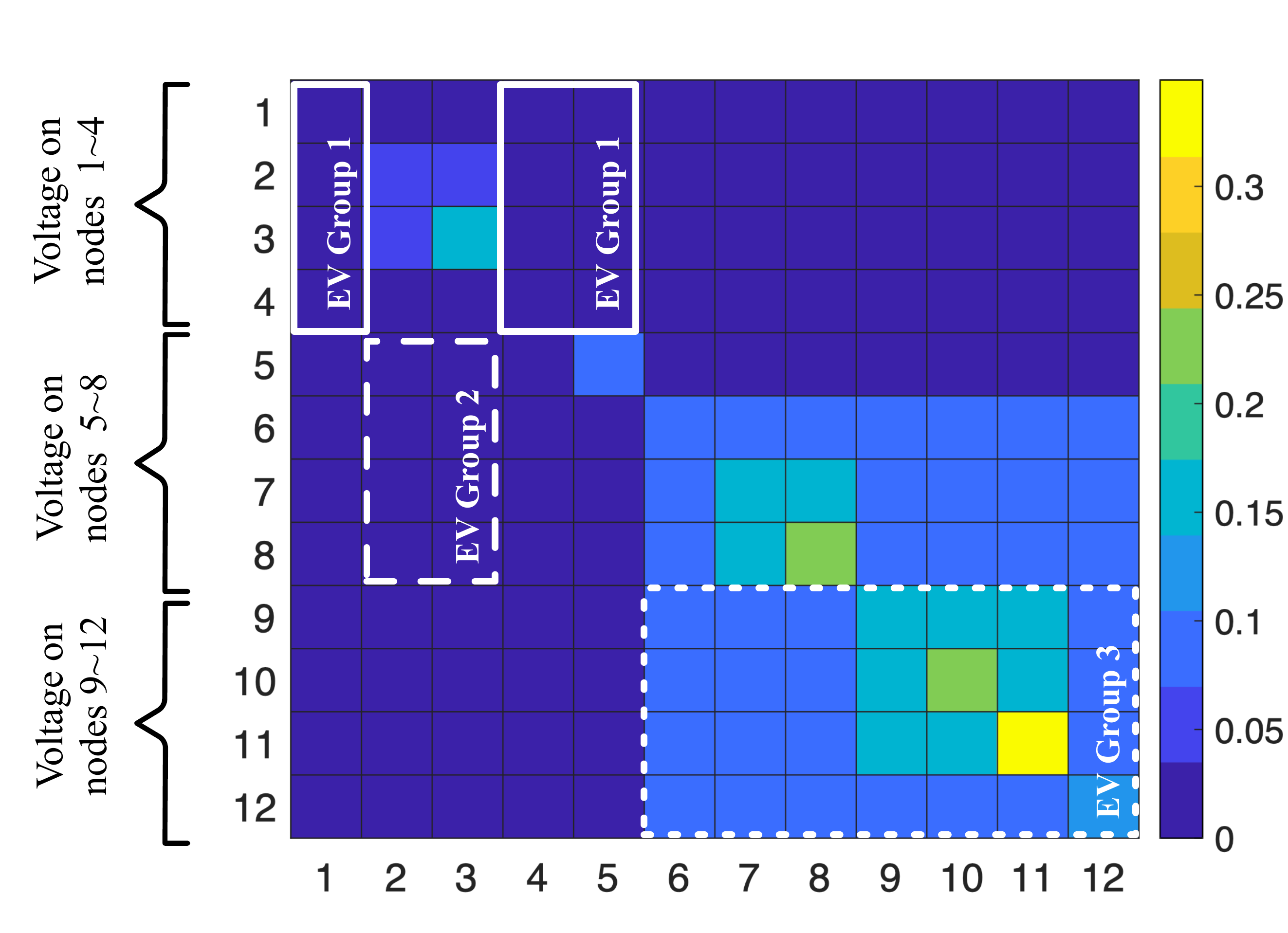}
    \caption{EV grouping and voltage subsets for the modified IEEE 13-bus test feeder (presented via the heatmap of the sensitivity matrix $\bm{R}$).}
    \label{grouping13}
\end{figure}

After grouping the EVs, voltage sub-vectors that reflect the major impacts should be constructed correspondingly. Here, we first take the modified IEEE 13-bus test feeder for example to demonstrate the voltage sub-vectors selection process, then provide the selecting rules. In Fig. \ref{grouping13}, the $\bm{R}$ matrix of a modified 13-bus test feeder is presented via heatmap. In total, we have 3 EV groups and 3 voltage subsets: EV Group 1 collects EVs connected at Nodes 1, 4, 5 and all EVs therein consider voltages at Nodes 1-4; Group 2 collects EVs connected at Nodes 2, 3 and all EVs therein consider voltages at Nodes 5-8; Group 3 collects EVs connected at Nodes 6-12 and all EVs therein consider voltages at Nodes 9-12. 

Let $d \in \mathbb{Z}_{0+}$ denote the dimension reduction  \emph{w.r.t.}  the original global voltage vector $\bm{V}(T)$, then the voltage sub-vector of the $s$th EV group can be written as
\begin{equation} 
    \bm{\hat{V}}_s(T) = \col\left(\hat{V}_{s,1}(T),~\cdots,~\hat{V}_{s,n-d}(T)\right)\in \mathbb{R}^{n-d}, 
    \label{4sn}
\end{equation}
for $s=1,\ldots,r$, where $r \in [1,n]$ is the total number of EV groups, and $\hat{V}_{s,l}(T)$, for $l=1,\ldots,n-d$, denotes the $l$th element of the voltage sub-vector of the $s$th EV group. Let $\tilde{\mathbb{V}}_{1},\ldots,\tilde{\mathbb{V}}_{r}$ denote the voltage subsets and $\tilde{\mathbb{V}}$ denote the full voltage set, then in order to make sure all $n$ nodal voltages are covered by in total $r$ EV groups, the following condition for grouping must be satisfied
\begin{equation}
\cup_{s=1}^{r}\tilde{\mathbb{V}}_s = \tilde{\mathbb{V}}.
\label{5sn}
\end{equation}
Moreover, to maximize the dimension reduction, the value of the dimension reduction $d$ can be determined by 
\begin{equation}
d = n - \lceil \frac{n}{r} \rceil, 
\label{6sn}  
\end{equation}
where $\lceil \cdot \rceil$ denotes the ceiling.

\textit{\textbf{Theorem 1:} \textbf{(Max Reduction and Least Overlapping)} Given the node number $n$ and the predefined group number $r$, by following \eqref{4sn}, \eqref{5sn} and \eqref{6sn}, the dimension reduction $d$ in \eqref{6sn} grants the maximum dimension reduction and minimum overlapping between voltage subsets. } \hfill $\blacksquare$

\emph{Proof:} The proof can be found in Appendix I. \hfill $\square$

In this particular case, 500 EVs connected to the modified IEEE 13-bus test feeder are divided into three groups, see Table \ref{table_example_13} for details.\begin{table}[!htb]
\caption{EV grouping methodology and voltage subset for each group in the modified IEEE 13-bus test feeder}
\label{table_example_13}
\begin{center}
\begin{tabular}{c|l|l|c}
\hline
EV Group Name & Location  & Voltage subset & Number of EVs\\
\hline
Group 1 & Nodes 1,4,5   & Nodes  1-4 & 100\\
Group 2 & Nodes 2,3   & Nodes   5-8 & 100\\
Group 3 & Nodes 6-12  & Nodes 9-12 & 300\\
\hline
\end{tabular}
\end{center}
\end{table}
We set the dimension reduction to $d=8$ to maximize the dimension reduction and make sure \eqref{5sn} is satisfied. The voltage sub-vectors of EV Groups 1-3 are
\begin{align*}
 \bm{\hat{V}}_1(T) &= \col\left(V_{1}(T),\cdots,V_{4}(T)\right),\\
 \bm{\hat{V}}_2(T) &= \col\left(V_{5}(T),\cdots,V_{8}(T)\right),\\
 \bm{\hat{V}}_3(T) &= \col\left(V_{9}(T),\cdots,V_{12}(T)\right).
\end{align*}

In this paper, we assume the voltage sub-vector of each EV group has the same dimension. The nonuniform partition will be discussed in our future work. Note that $\hat{\bm{V}}_1(T)$, $\hat{\bm{V}}_2(T)$, and $\hat{\bm{V}}_3(T)$ are not unique -- any pair of sub-vectors satisfying the grouping and dimension reduction rules \eqref{4sn}, \eqref{5sn} and \eqref{6sn} can be used. In real charging scenarios, each EV only considers the voltage subset of the group it belongs to so as to only consider its major impacts on the network, and overlapping between different voltage subsets is allowed but can be minimized when maximum dimension reduction is achieved.

\subsection{Shrunken Primal-Multi-Dual Subgradient Algorithm}

This section aims at developing a two-facet scalable decentralized algorithm to solve \eqref{eq:15} based on the dimension reduction technique developed in  Section \ref{section_dimension_reduction}. In our preliminary work \cite{2020cdcspmds}, we presented the outline of SPMDS. However, neither indication for the integration of dimension reduction into algorithm design nor the convergence analysis was provided in \cite{2020cdcspmds}. In this paper, we will largely enrich the discussion of the algorithm design and develop new theorems that provide convergence guarantees.

Through Lagrangian relaxation, the problem in \eqref{eq:15} can be solved by a primal-dual scheme in which EVs iteratively update their charging profiles. The relaxed Lagrangian of the problem in \eqref{eq:15} can be written as
\begin{equation}\label{Relaxed_Lag_rep} \mathcal{L}(\mathcal{U}, \bm{\lambda}) =\mathcal{F}(\mathcal{U})+\bm{\lambda}^{\mathsf{T}}\left( \bm{Y}_b - \sum_{i=1}^{v} \mathcal{D}_{i}\mathcal{U}_{i}\right),
\end{equation}
where $\bm{Y}_b = \underline{v}^{2} \bm{V}_{0} - \bm{V}_c$. Note that, local constraints in \eqref{eq:12sb} are not included in \eqref{Relaxed_Lag_rep}, but will be revisited in local updates.

Convexity of \eqref{eq:15} allows it to be solved through the fixed point problem \cite{boyd2004convex}
\begin{subequations}
\begin{align}
     \mathcal{U}_{i}^{*} &=\Pi_{\mathbb{U}_i}\left(\mathcal{U}_{i}^{*}-\nabla_{\mathcal{U}_{i}} \mathcal{L}\left(\mathcal{U}^{*},\bm{\lambda}^{*}\right)\right),\\
     \bm{\lambda}^{*} &=\Pi_{\mathbb{R}_{+}^{n K}}\left(\bm{\lambda}^{*}+\nabla_{\bm{\lambda}} \mathcal{L}\left(\mathcal{U}^{*}, \bm{\lambda}^{*}\right)\right),
\end{align}
\end{subequations}  
where $\Pi_{\mathbb{O}}(\bm{o})$ is the projection function projecting $\bm{o}$ onto the convex set $\mathbb{O}$, and in specific, 
\begin{subequations}\label{20}
\begin{align}
    \nabla_{\mathcal{U}_{i}} \mathcal{L}\left(\mathcal{U},\bm{\lambda}\right) &= \tilde{P}^{\mathsf{T}}(P_b+\tilde{P}\mathcal{U}_{i})+\rho \mathcal{U}_{i}+ \mathcal{D}_{i}^{\mathsf{T}}\bm{\lambda},\label{20a}\\
    \nabla_{\bm{\lambda}} \mathcal{L}\left(\mathcal{U}, \bm{\lambda}\right) &=\bm{Y}_b - \sum_{i=1}^{v} \mathcal{D}_{i}\mathcal{U}_{i}.\label{20b}
\end{align}
\end{subequations}

By following the proposed dimension reduction method and \textit{\textbf{Theorem 1}}, we group EVs into $r$ groups, and let $g_s$ denote the number of EVs in the $s$th group, for $s=1,\ldots,r$. By using the same $\hat{}$ notation as in \eqref{4sn} to denote the voltage sub-vector, then for the $s$th EV group, the reduced-dimension augmented slack node voltage vector is defined as $\hat{\bm{V}}_{0,s} = V_0^2\bm{1}_{(n-d)K}$. Correspondingly, the sub-vector $\bm{V}_c$ has its reduced-dimension counterpart $\hat{\bm{V}}_{c,s} \in \mathbb{R}^{(n-d)K}$. Further, let $\mathcal{D}_{d,i} \in \mathbb{R}^{(n-d) K \times K}$ and $\bm{D}_d \in \mathbb{R}^{(n-d) \times v}$  denote the reduced forms of $\mathcal{D}_{i}$ and $\bm{D}$, respectively. $\mathcal{D}_{d,i}$ and $\bm{D}_d$ are obtained by removing $d$ elements \emph{w.r.t.} the corresponding voltage sub-vectors from $\mathcal{D}_{i}$ and $\bm{D}$, respectively. Having the above definitions, we modify the subgradients in \eqref{20} to a reduced-dimension form as
\begin{subequations} \label{modified_gradient}
\begin{align}
    \tilde{\nabla}_{\mathcal{U}_{i}} \mathcal{L}\left(\mathcal{U},\bm{\lambda}_e\right) 
     &= \tilde{P}^{\mathsf{T}}(P_b+\tilde{P}\mathcal{U}_{i})+\rho \mathcal{U}_{i} + \mathcal{D}_{d,i}^{\mathsf{T}}\bm{\lambda}_e, \label{eq:20a}\\
    \tilde{\nabla}_{\bm{\lambda}_s} \mathcal{L}\left(\mathcal{U}, \bm{\lambda}_s\right)  
     &=\bm{\omega}_s \odot \bm{Y}_{b,s} - \sum_{\sum_{j=1}^{s-1}g_j+1}^{\sum_{j=1}^{s}g_j} \mathcal{D}_{d,i} \mathcal{U}_{i} = d_s(\mathcal{U}) \label{eq:20b}\\
    \bm{\lambda}_e &= \bm{\lambda}_1 + \cdots + \bm{\lambda}_r,\label{21c}
\end{align}
\end{subequations}
\noindent where $\odot$ denotes elementwise multiplication, $\bm{\lambda}_1,\cdots,\bm{\lambda}_r \in \mathbb{R}^{(n-d)K}$ denote the dual variables corresponding to the subsets of the global constraints of EV Groups $1$ to $r$, $\bm{\lambda}_e$ denotes the weighted dual variable in the modified primal subgradient \eqref{eq:20a}, $\bm{Y}_{b,s} \in \mathbb{R}^{(n-d)K}$ is a sub-vector of $\bm{Y}_{b}$ in \eqref{Relaxed_Lag_rep}, i.e., 
\begin{equation}
    \bm{Y}_{b,s} = \underline{v}^2\hat{\bm{V}}_{0,s}-\hat{\bm{V}}_{c,s},\nonumber
\end{equation}
and $\bm{\omega}_1,\bm{\omega}_2,\cdots,\bm{\omega}_r \in \mathbb{R}^{(n-d)K}$ denote the charging impact of each EV group on their designated voltage sub-vectors. 
The difference between \eqref{20a} and \eqref{eq:20a} lies in the last part, where $\mathcal{D}_{i}^{\mathsf{T}}\bm{\lambda}$ is replaced with $\mathcal{D}_{d,i}^{\mathsf{T}}\bm{\lambda}_e$. Also, we are creating multiple dual variables by modifying \eqref{20b} into \eqref{eq:20b}. Take EV Group 1 in the modified IEEE 13-bus test feeder as an example, the 100 EVs in Group 1 have impacts of \emph{magnitude} $\bm{\omega}_1$ on the voltage sub-vector $\bm{\hat{V}}_1 = \col( V_{1},\ldots,V_{7})$, i.e.,
global voltages of the node set \{$1,\ldots,7$\}. There may exist many ways to reflect this magnitude. Here, we define $\bm{\omega}_s \in [0,1]$, for $s=1,2,\ldots,r$ as
\begin{equation}\label{eq:22}
\bm{\omega}_s = \left(\sum_{\sum_{j=1}^{s-1}g_j+1}^{\sum_{j=1}^{s}g_j} \mathcal{D}_{d,i} \mathcal{U}_{i} \right) \oslash \left(\sum_{i=1}^{v} \mathcal{D}_{d,i} \mathcal{U}_{i}\right),
\end{equation}
where $\oslash$ denotes the elementwise division. In \eqref{eq:22}, the numerator denotes the charging influence of the EVs in the EV group $s$ on the $s$th voltage sub-vector while the denominator denotes the charging influence of all the EVs on the $s$th voltage sub-vector. In generic, \eqref{eq:22} can be interpreted as the portion of the impact of the agents in the group $s$ on the $s$th sub-vector.

\noindent \textbf{Remark 2:} A special case happens when $ \left[\sum_{i=1}^{v} \mathcal{D}_{d,i} \mathcal{U}_{i}\right]_{\jmath} = 0$ for some $\jmath$, where $[\cdot]_{\jmath}$ denotes the $\jmath$th entry. This indicates that no EV is charging at the $\jmath$th time slot. At this point, we define $\left[\bm{\omega}_1\right]_{\jmath} = \left[\bm{\omega}_2\right]_{\jmath} = \ldots = \left[\bm{\omega}_r\right]_{\jmath} = 1$ by directly using the upper bound of $\bm{\omega}_s$. Other options of $\bm{\omega}_s \in (0,1]$ may also be adopted without affecting the convergence. \hfill $\blacksquare$

\begin{figure}[!htbp]
    \centering
    \includegraphics[width=0.5\textwidth,trim = 0mm 20mm 5mm 5mm, clip]{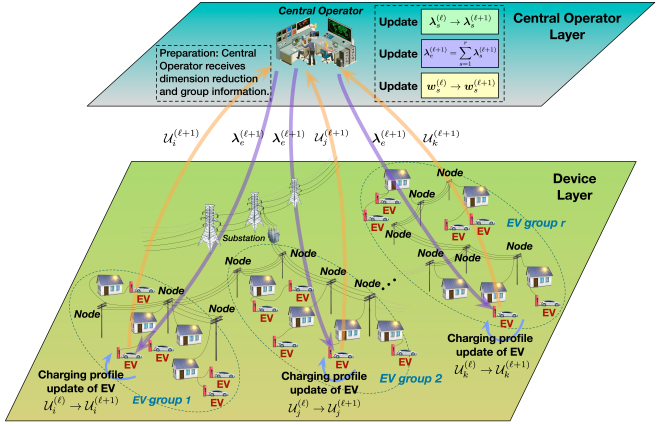}
    \caption{Information flow of decentralized EV charging control via SPMDS in a distribution network.}
    \label{information_flow}
\end{figure} 

Our previous work on SPDS \cite{liu2017decentralized} is dedicated to enabling a decentralized scheme for optimization problems with coupled objective functions and constraints, however its primal and dual updates suffer from large computational load caused by the dimension of the distribution network. In this paper, we propose a novel shrunken primal-multi-dual subgradient (SPMDS) algorithm based on the improved primal and dual gradients \eqref{modified_gradient} and with the integration of the developed dimension reduction technique to circumvent the drawbacks of SPDS. At the $\ell$th iteration, the proposed SPMDS updates the primal and multiple dual variables by following
\begin{subequations}\label{eq:23}
\begin{align}
\mathcal{U}_{i}^{(\ell+1)} &{=}\Pi_{\mathbb{U}_{i}}\left( \frac{1}{\tau_{\mathcal{U}}} 
\Pi_{\mathbb{U}_{i}}\left( \tau_{\mathcal{U}}\mathcal{U}_{i}^{(\ell)}-\alpha_{i, \ell} \tilde{\nabla}_{\mathcal{U}_{i}} \mathcal{L}\left(\mathcal{U}^{(\ell)}, \bm{\lambda}^{(\ell)}_{e}\right)\right)\right), \label{eq:23a}\\
\bm{\lambda}^{(\ell+1)}_s &{=} \Pi_{\mathbb{D}}\left( \frac{1}{\tau_{\mathcal{\lambda}_s}}\Pi_{\mathbb{D}}\left(
\tau_{\mathcal{\lambda}_s}
\bm{\lambda}^{(\ell)}_s+\beta_{s,\ell} \tilde{\nabla}_{\bm{\lambda}_s} \mathcal{L}\left(\mathcal{U}^{(\ell)}, \bm{\lambda}^{(\ell)}_s\right)\right)\right), \label{eq:23dual}\\
\bm{\lambda}^{(\ell)}_e &{=} \bm{\lambda}^{(\ell)}_1 + \cdots +  \bm{\lambda}^{(\ell)}_r, \label{eq:23b}
\end{align}
\end{subequations}
for $i=1,\ldots,v$ and $s=1,\ldots,r$, where $\bm{\lambda}^{(\ell)}_1,\bm{\lambda}^{(\ell)}_2, \ldots ,\bm{\lambda}^{(\ell)}_r$ denote the dual variables of EV Groups $1$ to $r$, respectively, $\mathbb{D}$ is the feasible set of $\bm{\lambda}^{(\ell)}_1,\bm{\lambda}^{(\ell)}_2, \ldots ,\bm{\lambda}^{(\ell)}_r$ (details can be found in \cite{liu2017decentralized}), and $\bm{\lambda}^{(\ell)}_e$ denotes the weighted dual variable of all $r$ groups which is universally adopted by all primal updates \eqref{eq:23a}. Adopting a universal dual variable $\bm{\lambda}^{(\ell)}_e$ in \eqref{eq:23a} has multifarious advantages: First, this can eliminate the necessity of CO transmitting particular dual variables to different EV groups; second, it assists the convergence of SPMDS; third, it reduces the computational cost. In addition, $\tau_{\mathcal{U}}\in [0,1]$ and $\tau_{\mathcal{\lambda}_1},\ldots,\tau_{\mathcal{\lambda}_r} \in [0,1]$ denote the shrinking parameters of the primal updates and $r$ dual updates, respectively. $\alpha_{i, \ell} > 0$, for $i=1,2,\ldots,v$, denotes the primal update step size of the $i$th EV, and 
$\beta_{s,\ell}>0$, for $s=1,2,\ldots,r$, denotes the dual update step size of the $s$th EV group. The convergence criterion is set correspondingly to the convergence of the charging profiles of all EVs and all the dual variables, i.e.,
\begin{equation}
    \epsilon=\left\|{\mathcal{U}}^{(\ell+1)}-\mathcal{U}^{(\ell)}\right\|_{2} + \sum_{s=1}^{r}
\left\Vert \bm{\lambda}_s^{(\ell+1)}-\bm{\lambda}_s^{(\ell)}\right\rVert_{2} < \epsilon_0,
    \label{eq:24}
\end{equation}
where $\epsilon_0$ is the predefined convergence tolerance. The complete SPMDS algorithm is summarized in Algorithm \ref{SPMDS_alrogithm_1}. Fig.  \ref{information_flow} presents the information flow of SPMDS for EV charging control in a distribution network described in \eqref{eq:23}. In the device layer, EVs are clustered into $r$ groups based on their locations in the distribution network, and each EV in a specific group is responsible for updating its own charging profile. While in the CO layer, the CO aggregates and controls the dual updates. During the primal and dual updates, CO receives the charging profiles $\mathcal{U}^{(\ell)}$ from all EVs and broadcasts the weighted dual variable $\bm{\lambda}_e^{(\ell)}$. This SPMDS-based framework only requires one CO for the multiple groups, leading to a lower communication and computation cost. Note that the CO needs to acknowledge the dimension reduction and group information before the iteration starts.

\begin{algorithm} 
\caption{SPMDS Algorithm}
\label{SPMDS_alrogithm_1}
\begin{algorithmic}[1]
\State Determine $r$ EV groups with $\mathcal{K}$-means using \eqref{3s} and the corresponding voltage subsets $\tilde{\mathbb{V}}_1,\ldots,\tilde{\mathbb{V}}_r$ by \eqref{4sn}, \eqref{5sn} and \eqref{6sn}.

\State Parameters initialization: CO initializes $\bm{\lambda}^{(0)}_{1},\ldots,$ $\bm{\lambda}^{(0)}_{r}$; EVs initialize $\mathcal{U}_i^{(0)}$; Tolerance $\epsilon_0$; Primal and dual step size of $\alpha_{i, \ell}>0$ and 
$\beta_{1,\ell}>0, \ldots, \beta_{r,\ell}>0$;  Shrinking parameters $\tau_{\mathcal{U}}$ and
$\tau_{\lambda_{1}},\ldots,\tau_{\lambda_{r}}$; Iteration counter $\ell=0$; Maximum iteration number $\ell_{max}$.

\While{$\epsilon > \epsilon_0$ and $\ell < \ell_{max}$}

\State The CO computes $\bm{\omega}_{1}^{(\ell)},\ldots,\bm{\omega}_{r}^{(\ell)}$ in \eqref{eq:22}, and calculates and broadcasts $\bm{\lambda}_{e}^{(\ell)}$ in \eqref{eq:23b} to all EVs.

\State Individual EV solves \eqref{eq:23a} and uploads its proposed charging schedule $\mathcal{U}_i^{(\ell+1)}$ to the CO.

\State The CO solves \eqref{eq:23dual} to obtain $\bm{\lambda}_{1}^{(\ell+1)},\ldots,\bm{\lambda}_{r}^{(\ell+1)}$.

\State  $\epsilon=\|{\mathcal{U}}^{(\ell+1)}-\mathcal{U}^{(\ell)}\|_{2}+ \sum_{s=1}^{r}
\| \bm{\lambda}_s^{(\ell+1)}-\bm{\lambda}_s^{(\ell)}\|_{2}$.

\State $\ell=\ell+1$.

\EndWhile
\end{algorithmic}
\end{algorithm}

\noindent {\bf{Remark 3:}} Though both SPMDS and SPDS emulate the public key encryption, where the dual variables $\bm{\lambda}_e$ and $\bm{\lambda}$ are the public keys while $\mathcal{D}_{d,i}$  and $\mathcal{D}_i$ are the private keys, SPMDS enhances the cyber-security of SPDS. In SPDS, once both $\bm{\lambda}$ and $\mathcal{D}_i$ are sniffed by cyber-attackers, they can be used to reverse engineer the distribution network topology and configuration. In contrast, this reverse engineering cannot be done in SPMDS with the reduced-dimension $\bm{\lambda}_e$ and the reduced-dimension $\mathcal{D}_{d,i}$. \hfill $\blacksquare$




\subsection{Convergence Analysis}

The theoretical foundation of SPMDS was not attempted in our preliminary work \cite{2020cdcspmds}. In this section, we will develop the optimality and convergence guarantees for SPMDS. Note that, the convexity of the problem under discussion guarantees the existence of the global optimum. The convergence analysis is grounded in Lyapunov stability theory. We will show that SPMDS has primal convergence, i.e., $\|\mathcal{U}^{(\ell+1)}-\mathcal{U}^{*}\|_{2}^{2} \to 0$ as $\ell \to \infty$, and dual convergence, i.e., $\sum_{s=1}^{r}\|\bm{\lambda}_s^{(\ell)}-\bm{\lambda}_s^{*} \|_{2}^{2} \to 0$ as $\ell \to \infty$.  We first present the convergence guarantees of SPMDS in \textit{\textbf{Theorem 2}} with uniform step sizes (i.e., $\alpha_{i,\ell}=\alpha, \beta_{s,\ell}=\beta$) for the clarity, then present  \textit{\textbf{Corollary 2.1}} to show the convergence in the nonuniform case.

\textit{\textbf{Theorem 2:}
Let $\{\zeta^{(\ell)}\}$ be a sequence generated by \eqref{eq:23}, where $\zeta^{(\ell)}$ is defined as 
\begin{equation}\label{11sss}
    \zeta^{(\ell)}\triangleq\col(\mathcal{U}^{(\ell)}, \bm{\lambda}_1^{(\ell)},\cdots, \bm{\lambda}_r^{(\ell)}).
\end{equation}
Let $\zeta^* = \col (\mathcal{U}^{*},\bm{\lambda}_1^*,\ldots,\bm{\lambda}_r^*)$ be the unknown optimizer, and define the Lyapunov function candidate
\begin{equation}\label{6s}
    \mathcal{V}{(\mathcal{U},\bm{\lambda}_1,\ldots,\bm{\lambda}_r)}\triangleq\beta^2 \left\Vert \mathcal{U}-\mathcal{U}^{*}\right\rVert_{2}^{2} + \alpha^2 \sum_{s=1}^{r}
\left\Vert \bm{\lambda}_s-\bm{\lambda}_s^{*}\right\rVert_{2}^{2},
\end{equation}
where $\mathcal{U},\bm{\lambda}_1,\ldots,\bm{\lambda}_r$ denote the algorithm states. If the parameters of SPMDS, i.e., $\alpha,\beta,\tau_\mathcal{U}, \tau_\lambda$ satisfy 
\begin{align}
    \max \left\{ \frac{M+\Psi L_\phi^2}{2 \Psi F_{\mathcal{U}}},\frac{N+ \Psi L_\phi^2}{2 \Psi F_{\lambda}}
    \right\} < \mu < 1,
    \label{8sss}
\end{align}
where 
\begin{align}
    M &= \frac{\alpha^2\beta^2}{\tau_{\mathcal{U}}^{2}}-\beta^2, ~N = \frac{\alpha^2\beta^2}{\tau_{\lambda}^{2}}-\alpha^2, \nonumber\\
    \Psi &= \max \left\{\frac{\alpha^2\beta^2}{\tau_{\mathcal{U}}^{2}},\frac{\alpha^2\beta^2}{\tau_{\lambda}^{2}}\right\}, \nonumber\\
    L_{\phi} &= \left\| \left[\rho+\left|\frac{\alpha-\tau_\mathcal{U}}{\alpha}\right|+L_{\nabla G}+rL_d,  \left|\frac{\beta-\tau_\lambda}{\beta}\right|+L_d \right] \right\|_2,\nonumber\\
    L_{\nabla G} &= nK \max_{i=1,\ldots,H}\{\bar{P}_i^2\},~L_{d} = HK\max _{j=1, \ldots, H K}\left\|D_{d,j}^{\mathsf{T}}\right\|_{2},\nonumber\\
    F_{\mathcal{U}} &= \rho+\frac{\alpha-\tau_\mathcal{U}}{\alpha}, ~F_{\lambda} = \frac{\beta-\tau_\lambda}{\beta},
\end{align}
then $\mathcal{V}{(\mathcal{U},\bm{\lambda}_1,\ldots,\bm{\lambda}_r)}$ is a Lyapunov function for the SPMDS and satisfies
\begin{align}
    &\mathcal{V}{(\mathcal{U}^{(\ell+1)},\bm{\lambda}_1^{(\ell+1)},\ldots,\bm{\lambda}_r^{(\ell+1)})} \leq \mathcal{V}{(\mathcal{U}^{(\ell)},\bm{\lambda}_1^{(\ell)},\ldots,\bm{\lambda}_r^{(\ell)})} \nonumber\\
    &+ A\left\|\mathcal{U}^{(\ell)}-\mathcal{U}^{*}\right\|_{2}^{2} +B\sum_{s=1}^{r}\left\Vert \bm{\lambda}_s^{(\ell)}-\bm{\lambda}_s^{*} \right\rVert_{2}^{2}
\end{align} 
where $A=M+\Psi L_\phi^2  -2\mu \Psi F_{\mathcal{U}}$, $B = N+ \Psi L_\phi^2 - 2\mu \Psi F_{\lambda}$, and $A<0,B<0$ are guaranteed. }  \hfill $\blacksquare$

\emph{Proof:} The proof can be found in Appendix II.  \hfill $\square$

\textit{\textbf{Corollary 2.1:}} \textit {\textbf{(Heterogeneous primal step sizes and dual step sizes)} SPMDS allows the agents to independently choose their own primal step sizes, and the CO to choose independent dual step size for each virtual EV group. Let the $i$th agent select its own primal step size $\alpha_{i,\ell}$, and the $s$th EV group has its dual step size $\beta_{s,\ell}$, then SPMDS achieves convergence if there exits $\alpha_{l,\ell}$ and $\beta_{l,\ell}$ that minimize the right hand side of \eqref{26ss} and $\alpha_{u,\ell}$ and $\beta_{u,\ell}$ that maximize the right hand side of \eqref{41ss}.} \hfill $\blacksquare$

\emph{Proof:} Let $\mathbb{P}$ and $\mathbb{U}$ denote the primal step size set and dual step size set, respectively. Thus, there always exist such $\alpha_{l,\ell} \in \mathbb{P}$ and $\beta_{l,\ell} \in \mathbb{U}$ that minimize the right hand side of \eqref{26ss}, and $\alpha_{u,\ell} \in \mathbb{P}$ and $\beta_{u,\ell} \in \mathbb{U}$ that maximize the right hand side of \eqref{41ss}. Then, by following the same proof procedure for \textit{\textbf{Theorem 2}}, \textit{\textbf{Corollary 2.1}} can be proved.   \hfill $\square$

\noindent {\bf{Remark 4:}} Though the NP-hard nature of $\mathcal{K}$-means offers a local optimum of the clustering, the convergence of SPMDS is guaranteed under any sub-optimal solutions obtained from $\mathcal{K}$-means. Besides, the least overlapping in \textbf{\textit{Theorem 1}} is not a necessary condition for the convergence of SPMDS, namely as long as \eqref{4sn}, \eqref{5sn} and \eqref{6sn} are satisfied, overlapping in voltage sub-vectors is 
optional. \hfill $\blacksquare$

\noindent {\bf{Remark 5:}} \textit{\textbf{(Parameter selection guidance)}} \textbf{\textit{Theorem 2}} and \textbf{\textit{Corollary 2.1}} give the rules of convergence for SPMDS under a proper set of parameters, i.e.,  $\alpha$ ($\alpha_{i,\ell}$), $\beta$ ($\beta_{s,\ell}$), $\tau_\mathcal{U}$, and $\tau_\lambda$. For the easiness of parameter tuning, one may select sufficiently small primal update step size $\alpha$ ($\alpha_{i,\ell}$) and dual update step size $\beta$ ($\beta_{s,\ell}$) to guarantee the convergence criteria in \textbf{\textit{Theorem 2}} or \textbf{\textit{Corollary 2.1}} first, then gradually increase them for a faster convergence speed till divergence. \hfill $\blacksquare$

\subsection{Computational Load Analysis}

SPMDS and SPDS have the similar primal-dual structure, but the former requires less memory and less computational cost for each EV charger (or equivalent onboard controller), as well as less computation time for the CO. In what follows, we analyze the computational cost reduction enabled by SPMDS. 

The computation cost difference between \eqref{20a} and  \eqref{eq:20a} lies in the last term $\mathcal{D}_{i}^{\mathsf{T}}\bm{\lambda}$ in \eqref{20a}. Without the EV grouping and dimension reduction strategy, $\mathcal{D}_{i}^{\mathsf{T}}\bm{\lambda}$ requires $(2nK-1)K$ floating point operations (FLOPS) from each EV; after the dimension reduction of $d$ for each group, the FLOPS reduce to $2(n-d)K^2-K$. Both \eqref{20a} and \eqref{eq:20a} have the same FLOPS of $4K^2-K$ for the first term $\tilde{P}^{\mathsf{T}}(P_b+\tilde{P}\mathcal{U}_{i})$. Therefore the total FLOPS reduction and the corresponding reduction ratio of the primal gradient  calculation in \eqref{eq:20a} are 
\begin{subequations}\label{eq:25sn}
\begin{align}
    \mathcal{F}_{pt} &= (2nK-1)K - ((2n-d)K-1)K  \nonumber\\
    &= 2dK^2,\label{eq:25a}\\
    \mathcal{F}_{pr} &= \frac{2dK^2}{(2nK-1)K+4K^2-K} \nonumber \\
    &= \frac{dK}{(n+2)K-1}.\label{eq:25b}
    \end{align}
\end{subequations}
Eqn. \eqref{eq:25a} shows a proportional relation between the absolute computational cost reduction and the dimension reduction $d$.
\begin{figure*}[!htb]
  \centering
  \begin{tabular}[b]{c}
    \includegraphics[width=.45\linewidth,height =.25\linewidth]{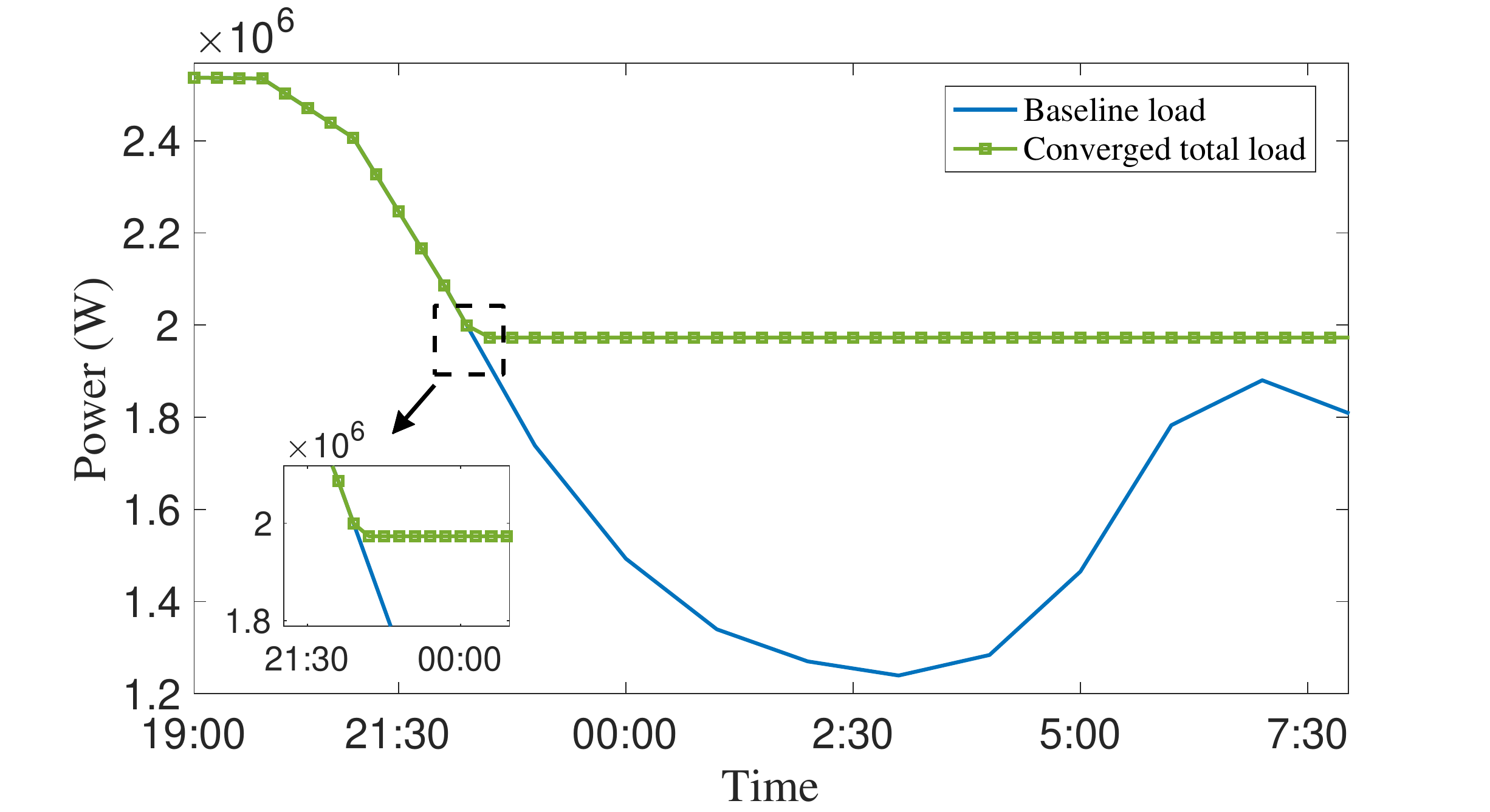} \\
    \small (a)
  \end{tabular} \qquad
  \begin{tabular}[b]{c}
    \includegraphics[width=.45\linewidth,height =.25\linewidth]{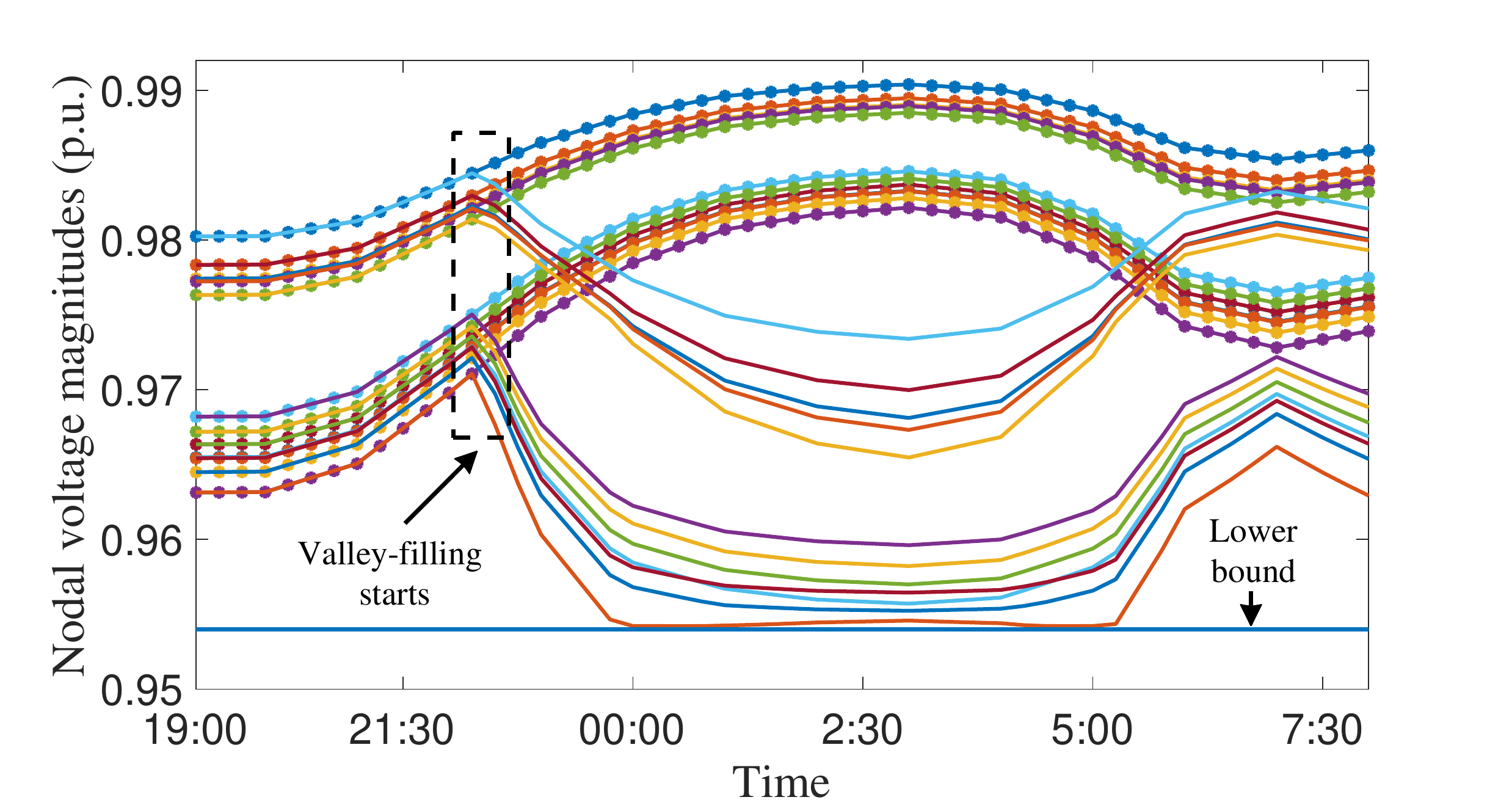} \\
    \small (b)
  \end{tabular}
    \begin{tabular}[b]{c}
    \includegraphics[width=.45\linewidth,height =.25\linewidth]{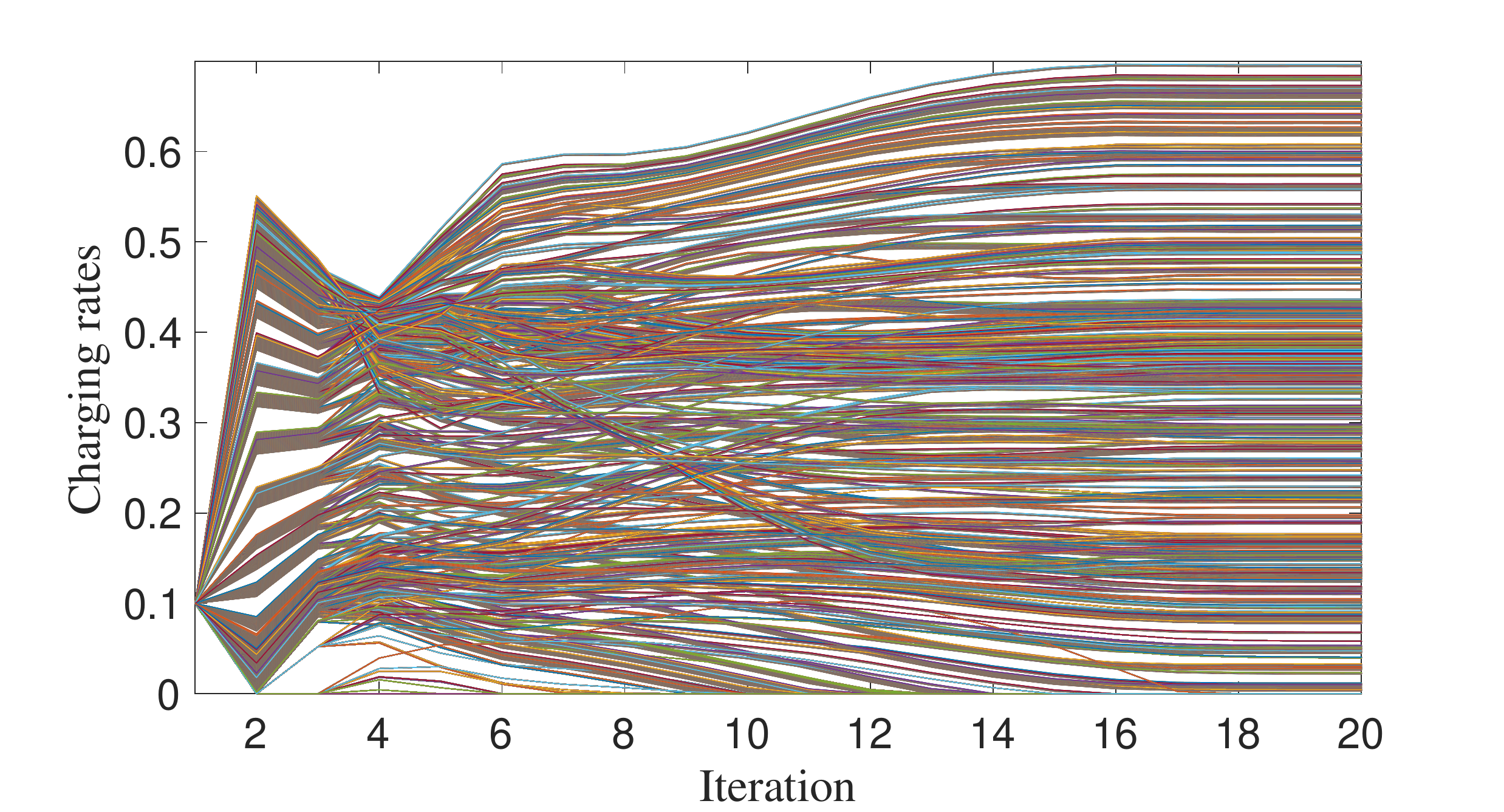} \\
    \small (c)
  \end{tabular} \qquad
  \begin{tabular}[b]{c}
    \includegraphics[width=.45\linewidth,height =.25\linewidth]{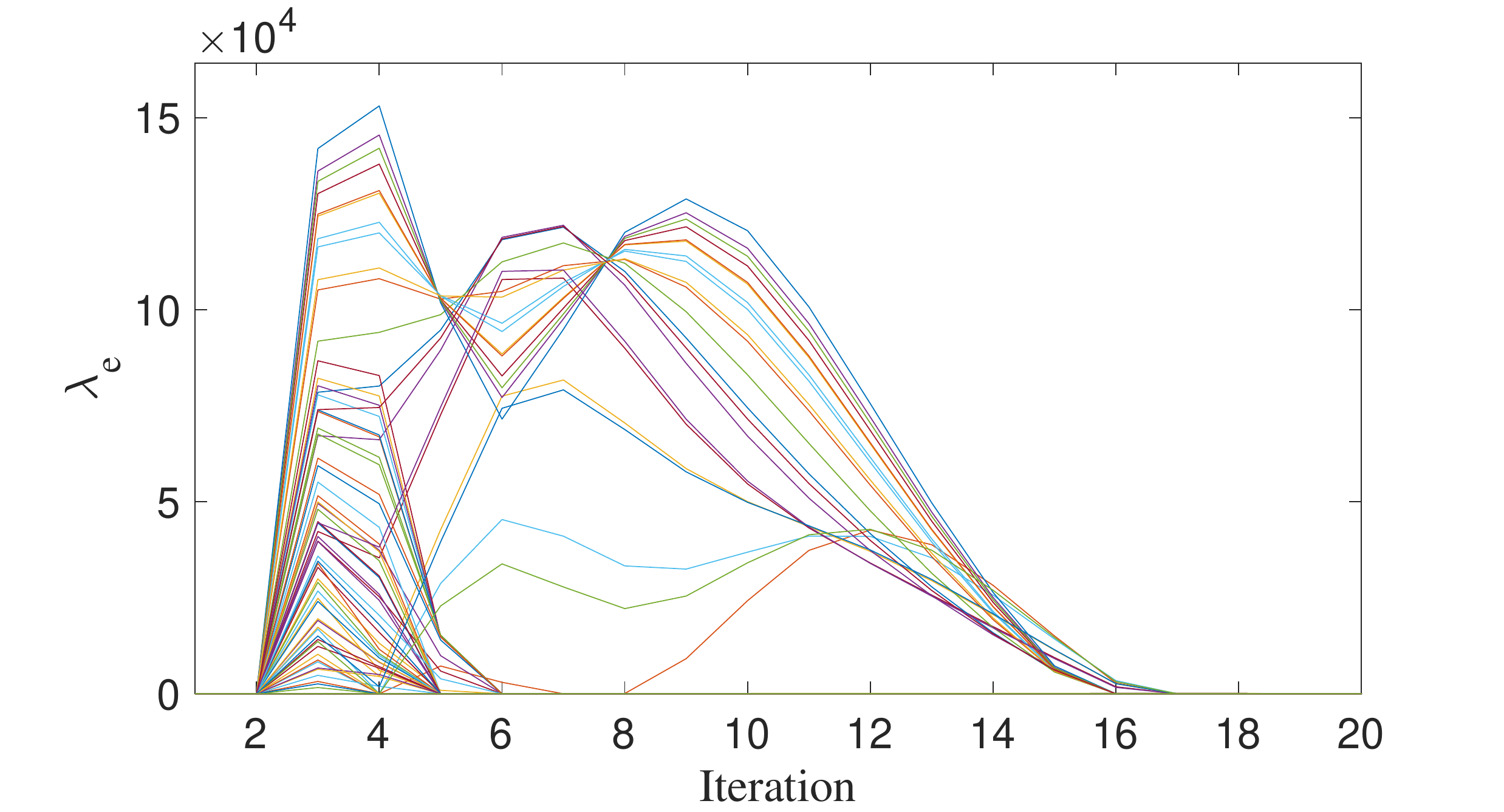} \\
    \small (d)
  \end{tabular}
  \caption{Valley-filling results for 500 EVs on the modified IEEE 13-bus test feeder (20 iterations) (a) Baseline load (solid line) and total load at the 20th iterations (dashed line) (b) Nodal voltage magnitudes of baseline load (dashed lines) and total load (solid lines) (c) Charging rates evolution of the 500 EVs in 20 iterations (d) Convergence of the weighted dual variable $\bm{\lambda}_e$ in 20 iterations }
\label{sim_IEEE_13}
\end{figure*}

\noindent For the dual update \eqref{eq:23b}, $\bm{\lambda}_1,\cdots,\bm{\lambda}_r$ can be calculated in parallel. Therefore, the total calculation time for the dual update depends on (i) the time consumed by the EV group that has the largest number of EVs, where we assume this group to be  $r_{m}$ with $g_{m}$ EVs, and (ii) the calculation time for \eqref{eq:20b}. The FLOPS in the dual subgradient calculation in \eqref{20b} is $2vnK^2$, and the FLOPS in the dual subgradient calculation in \eqref{eq:20b} is $(2v_mK+1)(n-d)K$. Note that the CO needs extra $(r-1)(n-d)K$ FLOPS due to \eqref{21c}. However, this is relatively small and negligible compared to the FLOPS of either \eqref{20b} or \eqref{eq:20b}. Therefore, when analyzing the computational cost, we can ignore \eqref{21c} and only consider the computational cost of Group $r_{m}$. Consequently, the cost reduction for the  dual gradient calculation in \eqref{eq:20b}  \emph{w.r.t.} FLOPS are
\begin{subequations}
\begin{align}
    \mathcal{F}_{dt} &= 2vnK^2 - (2v_mK^2+K)(n-d)\nonumber\\
    &= (2vn-2v_m(n-d))K^2 - K(n-d),
    \label{eq:26a}\\
    \mathcal{F}_{dr} &= \frac{(2vn-2v_m(n-d))K^2 - K(n-d)}{2vnK^2}.
    \label{eq:26b}
\end{align}
\label{eq:26sn}
\end{subequations}

\noindent \textbf{Remark 6:} For a larger-scale distribution network and a longer time interval $K$ (unnecessarily restricted to valley-filling problem), the dimension reduction $d$ plays a critical role in computational cost reduction. Besides, the dimension reduction makes it possible to decrease the memory cost, therefore taking advantage of the micro controller units embedded in the EV chargers. This will largely facilitate the deployment of the decentralized EV charging control architecture. By using \textit{\textbf{Theorem 1}}, we can have the maximum reduction which leads to the maximum reduction of FLOPS in the primal update as $\mathcal{F}_{pt} = 2(n - \lceil \frac{n}{r}\rceil)K^2$ and $ \mathcal{F}_{pr} = (n - \lceil \frac{n}{r} \rceil)K/((n+2)K-1)$, and the maximum reduction of FLOPS in the dual udpate as $\mathcal{F}_{dt} = (2vn-2v_m\lceil \frac{n}{r} \rceil)K^2 - K\lceil \frac{n}{r} \rceil$ and $ \mathcal{F}_{dr} = ((2vn-2v_m\lceil \frac{n}{r} \rceil)K^2 - K\lceil \frac{n}{r} \rceil)/(2vnK^2)$. \hfill $\blacksquare$

\section{Simulation Results}

In this section, we first conduct simulations of EV charging control on a modified IEEE 13-bus test feeder and a modified IEEE 123-bus test feeder, then test with a traffic congestion optimization problem to further manifest the generality of the proposed approaches.

\subsection{Decentralized EV Charging Control}
We consider two scenarios to verify the efficacy and efficiency of the proposed SPMDS-based decentralized EV charging control for valley-filling. First, we introduce some common parameters for both scenarios. The valley-filling period starts at 19:00 and ends at 8:00 next day, with 15-min time intervals and $K=52$ time slots correspondingly. The baseline load is assumed to be known and is directly adopted from \cite{liu2017decentralized}. All EVs have the maximum charging power of $6.6$ kW, EVs' random charging requirements vary from $20\%$ to $60\%$ state of charge, and the uniform charging efficiency is $\eta_i = 0.9$. The lower bound of the voltage is set to be $\underline{v} = 0.954$, which is slightly higher than the ANSI C84.1 standard to compensate for the LinDistFlow model inaccuracy. The slack node voltage magnitude is $V_0 = 4.16$ kV for both scenarios. 
\subsubsection{Scenario 1: Modified IEEE 13-bus Test Feeder} The modified IEEE 13-bus test feeder is adopted from our previous work \cite{liu2017decentralized}.
In this case, 50 EVs are connected at each node except that Nodes 1 and 6 have no load. The grouping and voltage subset selection results have been previously shown in Table \ref{table_example_13}. The primal step sizes and the dual step sizes are set uniformly as $\alpha = 2.8 \times 10^{-10}$  and $\beta = 1.8$. The shrinking parameters are empirically chosen to be $\tau_{\mathcal{U}} =\tau_{\bm{{\lambda}}_1} =\tau_{\bm{\lambda}_2}=\tau_{\bm{\lambda}_3} = 0.98$. 

The simulation results of 20 iterations are presented in Fig. \ref{sim_IEEE_13}. Fig. \ref{sim_IEEE_13}(a) shows the baseline load in contrast to the total load with the participation of 500 EVs -- the valley caused by the overnight baseline load reduction is filled by the controlled EV charging load between 22:00 and 8:00 the next day. Fig. \ref{sim_IEEE_13}(b) shows that the voltage magnitudes of all nodes are well maintained above the lower bound $0.954$ p.u. Fig. \ref{sim_IEEE_13}(c) presents the primal convergence of all EVs' charging rates and Fig. \ref{sim_IEEE_13}(d) shows the dual convergence of the weighted dual variable $\bm{\lambda}_e$.

The dimension reduction is set to the maximum $d = 8$ according to \textit{\textbf{Theorem 1}}. As a result, comparing with the SPDS, $\mathcal{F}_{pt} = 43,264$ FLOPS are reduced in each primal update, i.e., $\mathcal{F}_{pr} = 57.22 \%$. For each dual update, EV Group 3 has the largest number of EVs, i.e., $g_m = 300$. Then the FLOPS reduction in the dual subgradient calculation is $\mathcal{F}_{dt} = 25,958,192$, i.e., $\mathcal{F}_{dr} = 80\%$ for one iteration.

\subsubsection{Scenario 2: Modified IEEE 123-bus Test Feeder}
To verify the scalability of the proposed SPMDS algorithm, a 4.16kV modified IEEE 123-bus test feeder is considered \cite{IEEEPES}, where the voltage regulators are removed from the original IEEE 123-bus test feeder \cite{hu2020voltage} to validate the effectiveness of the proposed algorithm, and the impedance of segments are taken from phase A. In particular, we renumbered the network and broke the connection between Node 38 and Node 116 in the new numbering system to make it a radial distribution network. In addition, 5 EVs are connected at each node except that Nodes 1 and 6 have no load.

\begin{figure}[!htbp]
    \centering
    \includegraphics[width=0.45\textwidth,trim = 0mm 0mm 0mm 0mm, clip]{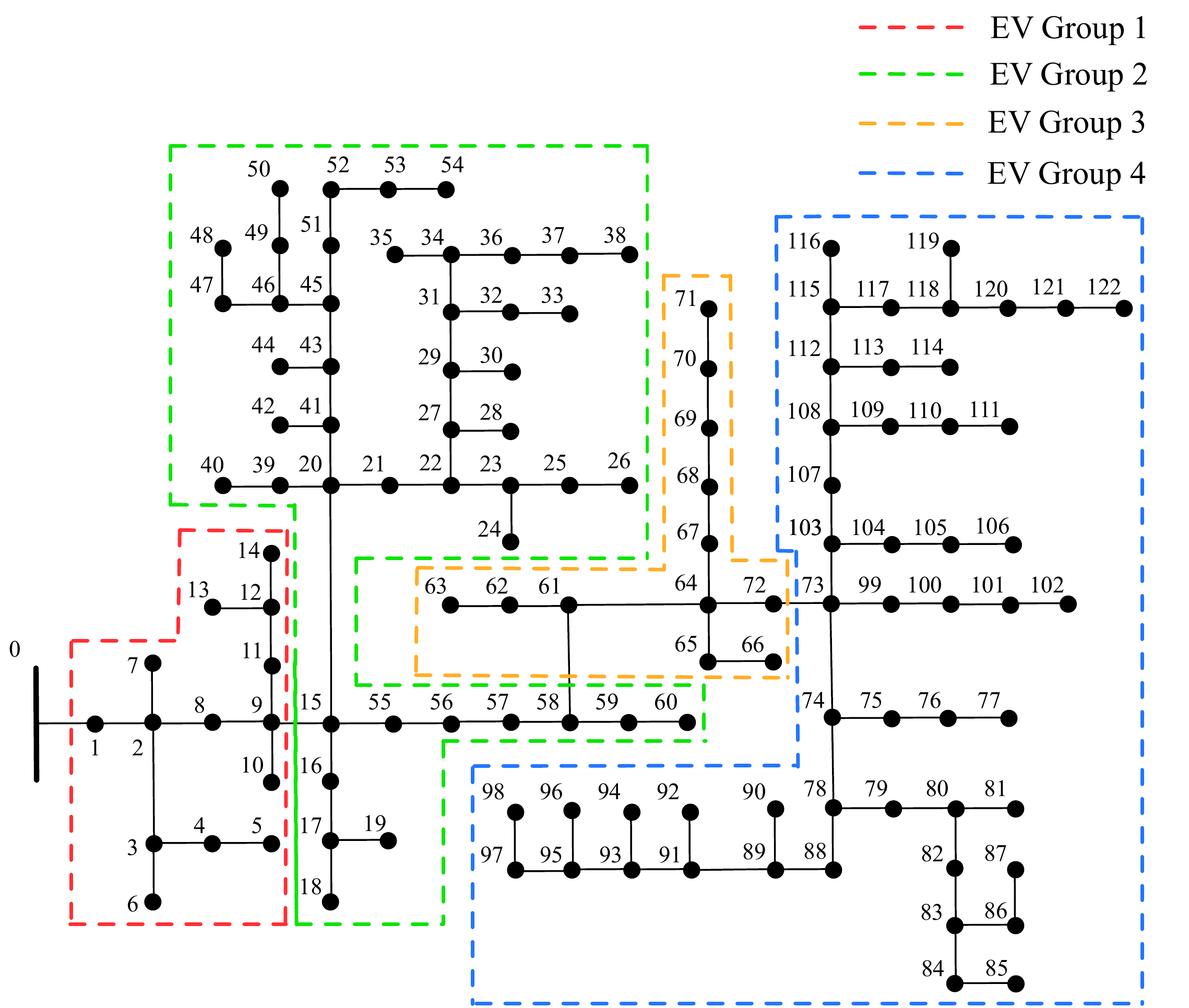}
    \caption{Network partitioning strategy via $\mathcal{K}$-means of the modified IEEE 123-bus test feeder (with clustering  number $r = 4$).}
    \label{sim_IEEE_123_distribution}
\end{figure} 

\begin{figure}[!htbp]
    \centering
    \includegraphics[width=0.4\textwidth,trim = 0mm 0mm 0mm 0mm, clip]{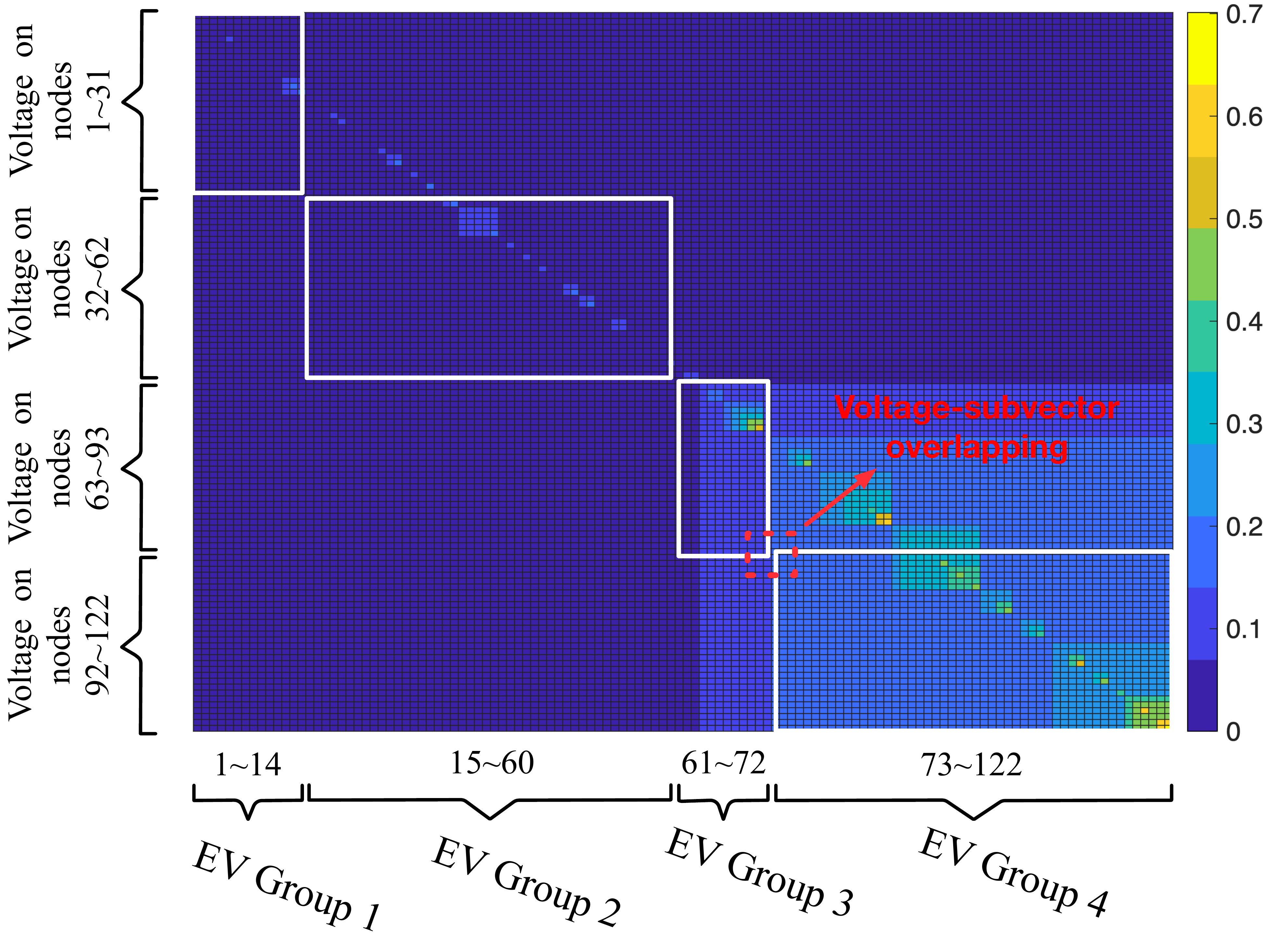}
    \caption{EV grouping and voltage subsets for the modified IEEE 123-bus test feeder (presented via the heatmap of the sensitivity matrix $\bm{R}$).}
    \label{sim_IEEE_123_heat}
\end{figure}

Fig. \ref{sim_IEEE_123_distribution} presents that total 122 nodes are partitioned into 4 clusters, and 600 EVs connected at the corresponding nodes can be divided into 4 groups accordingly. The EV grouping and voltage subset selections are presented in Fig. \ref{sim_IEEE_123_heat} via the heatmap of the sensitivity matrix $\bm{R}$, note that an overlapping of the voltage subsets between the EV group 3 and EV group 4 exists, i.e., Nodes 92 and 93 are in both $\tilde{\mathbb{V}}_{3}$ and $\tilde{\mathbb{V}}_{4}$. The dimension reduction is set to $d=91$ with the maximum dimension reduction by following \textit{\textbf{Theorem 1}}.  Table \ref{table_example_123} \begin{table}[!htb]
\caption{EV grouping methodology and voltage subset for each group in the modified IEEE 123-bus test feeder}
\label{table_example_123}
\begin{center}
\begin{tabular}{c|l|l|c}
\hline
EV Group Name & Location & Voltage subset & Number of EVs\\
\hline
Group 1 & Nodes 1-14    & Nodes  1-31 & 60\\
Group 2 & Nodes 15-60   & Nodes  32-62 & 230\\
Group 3 & Nodes 61-72   & Nodes  63-93& 60\\
Group 4 & Nodes 73-122  & Nodes 92-122 &250\\
\hline
\end{tabular}
\end{center}
\end{table}presents the details of the EV grouping and voltage subsets. The primal step sizes and the dual step sizes are set uniformly as $\alpha = 3 \times 10^{-10}$ and $\beta = 0.1$. The shrinking parameters are $\tau_{\mathcal{U}} =0.98, \tau_{\bm{{\lambda}}_1}=\tau_{\bm{{\lambda}}_2}=\tau_{\bm{{\lambda}}_3} =\tau_{\bm{{\lambda}}_4} = 0.97$. 

\begin{figure*}[!htbp]
  \centering
  \begin{minipage}{.45\textwidth}
    \centering
    \includegraphics[width=\linewidth]{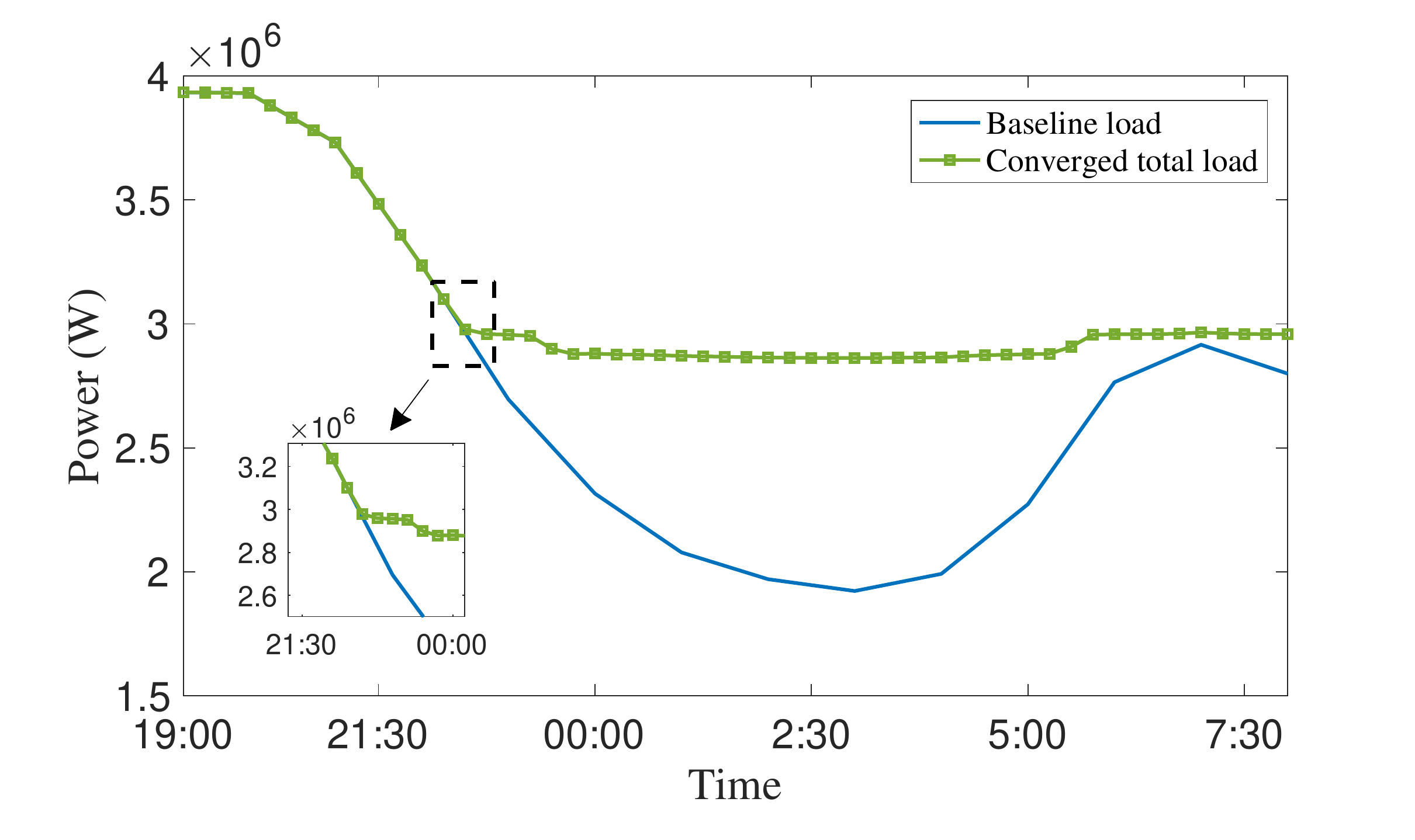}
    \small (a)
  \end{minipage} \quad
  \begin{minipage}{.45\textwidth}
    \centering
    \includegraphics[width=\linewidth]{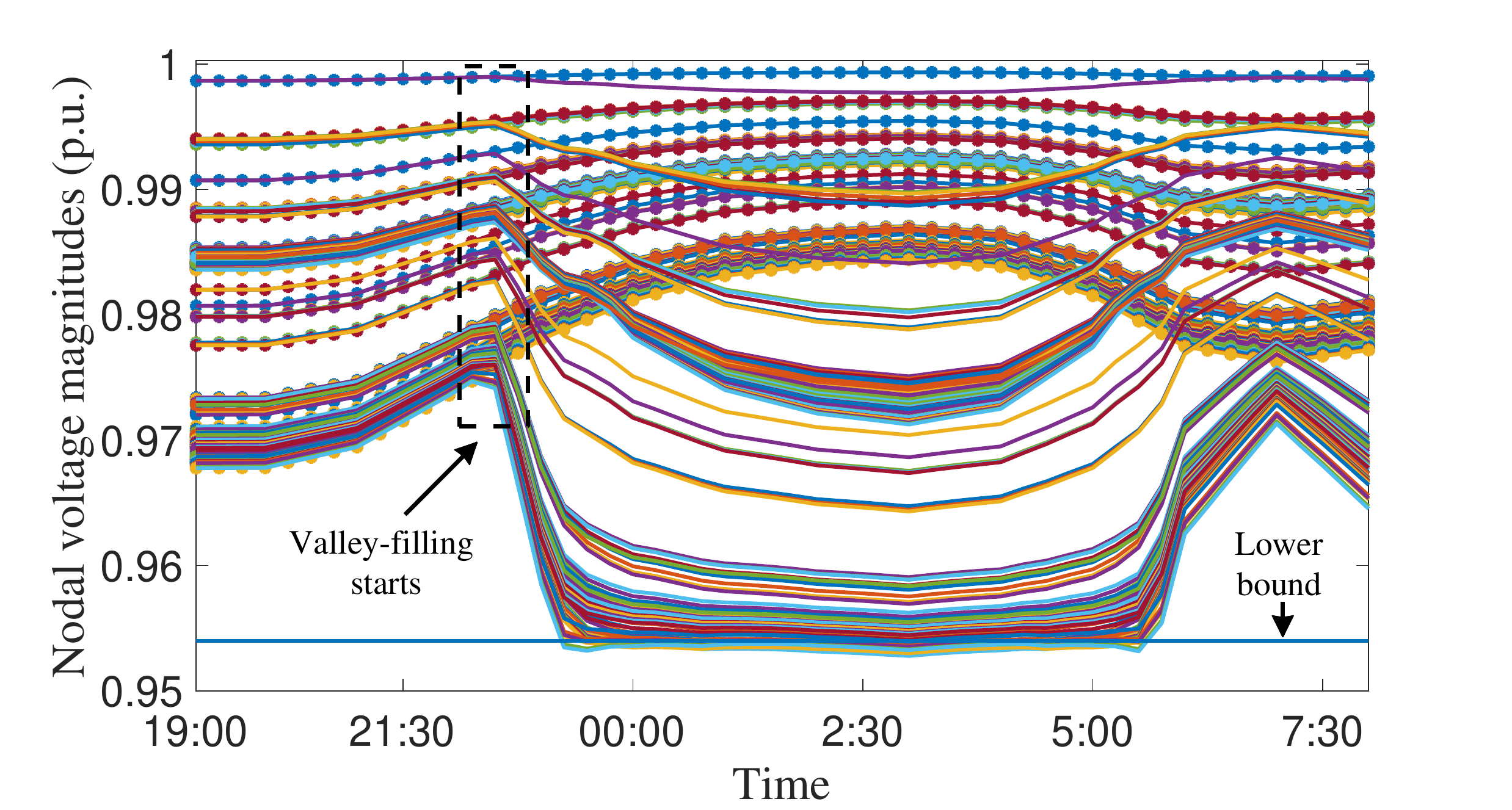}
    \small (b)
  \end{minipage} \quad
    \caption{Valley-filling results for 600 EVs on the modified IEEE 123-bus test feeder (30 iterations) (a) Baseline load (solid line) and total load at the 30th iterations (dashed line) (b) Nodal voltage magnitudes of baseline load (dashed lines) and total load (solid lines)}
  \label{sim_IEEE_123}
\end{figure*}

Fig. \ref{sim_IEEE_123}(a) and Fig. \ref{sim_IEEE_123}(b) depict the valley-filling performance and the voltage control, respectively, which clearly indicate the efficacy of the proposed method. After 30 iterations, though Fig. \ref{sim_IEEE_123}(b) exhibits subtle voltage violations due to numerical calculations, the converged results are already good enough to satisfy the engineering use. The engineering remedy to the subtle violation could be tightening the constraint bounds, which was already integrated by setting a slightly higher voltage lower bound.


Because of the dimension reduction $d=91$, $\mathcal{F}_{pt} = 492,128$ FLOPS are reduced in each primal update, i.e. $\mathcal{F}_{pr} = 73.4 \%$ FLOPS reduction compared to the full dimensional case. For the dual update, EV Group 4 has the largest number of EVs, i.e., $g_m = 250$, and $\mathcal{F}_{dt} = 353,951,988$ FLOPS in each dual update are reduced, i.e.,  $\mathcal{F}_{dr} = 89.41\%$ FLOPS reduction for one iteration compared to the full-dimension case.

\subsection{Transportation Congestion Control} 

To better illustrate the applicability and generality of SPMDS, we consider a transportation congestion optimization problem over a network with $N$ agents and $L$ links. Suppose each agent $i$ travels along a route with transmission rate $x_i$, then the shared paths between all agents arise traffic congestion. The goal is to optimize the congestion and minimize the congestion cost, which is an NOP.

\begin{figure}[!htbp]
    \centering
    \includegraphics[width=0.26\textwidth,trim = 5mm 5mm 5mm 5mm, clip]{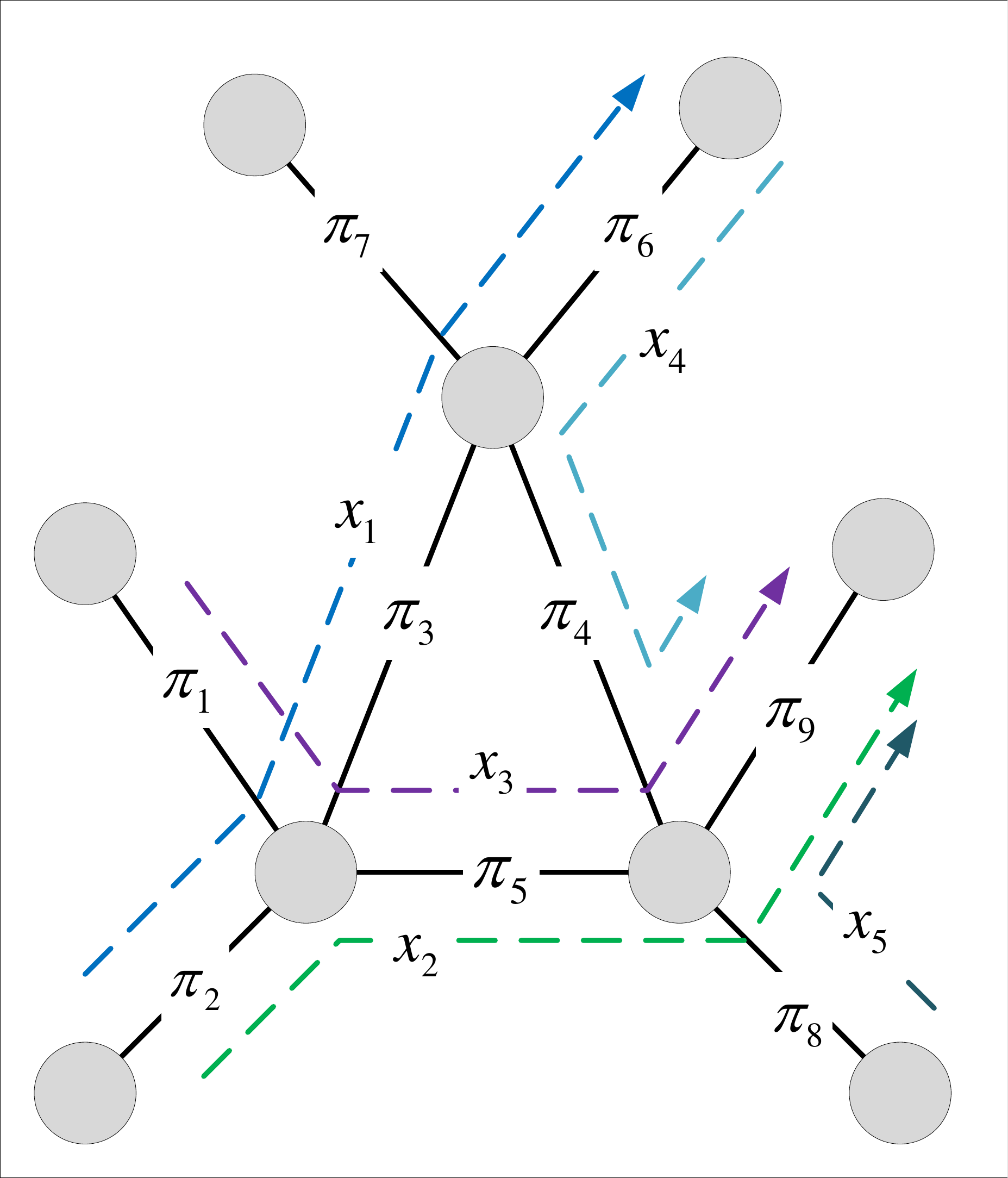}
    \caption{A transportation network with 9 links and 5 agents \cite{koshal2011multiuser}.}
    \label{traffic_flow}
\end{figure}

An example of a transportation network with 5 users and 9 links is shown in Fig. \ref{traffic_flow}. This example was used in \cite{koshal2011multiuser} to demonstrate RPDS.
Let $\pi_j$ denote the $j$th link and $\bm{x} = \col(x_1,\dots,x_N)$. The aggregated user utility needs to be maximized by varying $\bm{x}$ \cite{he2007towards}, and the utility function of agent $i$ is given by \cite{kelly1998rate}
\begin{equation}
    f_i(x_i) = k_i\log (1+x_i),
    \label{32sn}
\end{equation}
then the coupled congestion cost arsing from same usage of links across all agents is defined as 
\begin{equation}
    c(\bm{x}) = \sum_{i=1}^{N} \sum_{l \in L}x_{li}\sum_{m=1}^{N}x_{lm},
    \label{33ss}
\end{equation}
where $x_{lj}$ is the flow of agent $j$ on link $l$. 
Therefore, the objective function is coupled by the congestion cost across all agents in  \eqref{33ss}. Consequently, the transportation congestion optimization problem can be formulated as 
\begin{equation}
\begin{aligned}
& \underset{\bm{x}}{\text{min}} & & {\sum_{i=1}^{N} -f_i(x_i) + c(\bm{x})} \\
& \text{s.t.} & &  \bm{x} \geq 0\\
& & &  \bm{A}\bm{x} \leq \bm{b},
\label{32ssss}
\end{aligned}
\end{equation} 
where $\sum_{i=1}^{N}-f_i(x_i)$ comes from the maximization of agents' utility functions over transmission rates $\bm{x}$, constraint $\bm{A}\bm{x} \leq \bm{b}$ represents agent traffic rates over the network, $\bm{A} \in \mathbb{R}^{L\times N}$ denotes the link-route incidence matrix, i.e., $\bm{A}_{ji}=1$ if the link $j$ is on the path of agent $i$, and $\bm{A}_{ji}=0$ otherwise, and $\bm{b}$ is the link capacity vector with $b_j$ denoting the maximum aggregate traffic through link $l$.
\begin{table}[!htb]
\caption{Traffic network and agent data}
\label{table_example_congestion_flow}
\begin{center}
\begin{tabular}{c|c|c}
\hline
Agent Name & Links traversed & $k_i$\\
\hline
 1 & $\lambda_2,\lambda_3,\lambda_6$    & 10 \\
 2 & $\lambda_2,\lambda_5,\lambda_9$   & 0 \\
 3 & $\lambda_1,\lambda_5,\lambda_9$   & 10\\
 4 & $\lambda_6,\lambda_4,\lambda_9$  & 10 \\
 5 & $\lambda_8,\lambda_9$  & 10 \\
\hline
\end{tabular}
\end{center}
\end{table}

Table \ref{table_example_congestion_flow} concludes the traffic flow in the network and $k_i$ of the utility function \eqref{32sn}, the link capacity vector is set as $\bm{b} = \bm{1}_9$. 
The  details  of  the  grouping  and  congestion  subsets are shown in Table \ref{table_example_congestion_group}. 

\begin{table}[!htb]
\caption{Grouping methodology and congestion subset for each group in the transportation network}
\label{table_example_congestion_group}
\begin{center}
\begin{tabular}{c|c|c}
\hline
Group Name & Group members & Congestion sub-vector\\
\hline
 Group 1 & Agents 1,2   & $\bm{b}_1 = [1,1,1,1,1]^{\mathsf{T}}$ \\
 Group 2 & Agents 3-5   & $\bm{b}_2 = [1,1,1,1,1]^{\mathsf{T}}$  \\
\hline
\end{tabular}
\end{center}
\end{table}

\begin{figure}[!tbp]
    \centering
    \includegraphics[width=0.45\textwidth,trim = 0mm 0mm 0mm 0mm, clip]{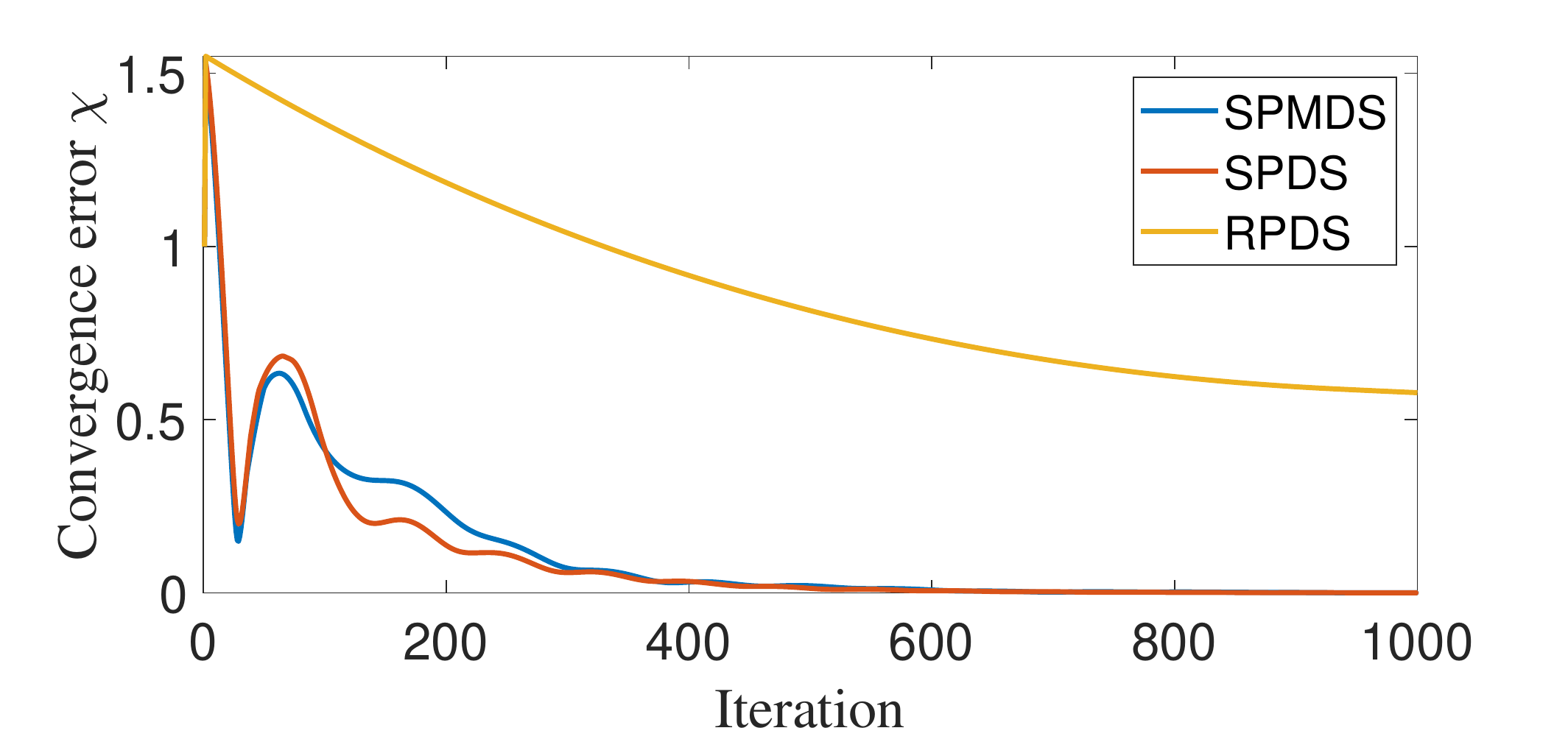}
    \caption{Convergence errors comparison between SPMDS, SPDS and RPDS.}
    \label{traffic_compare}
\end{figure}

\begin{figure}[!tbp]
    \centering
    \includegraphics[width=0.45\textwidth,trim = 0mm 0mm 0mm 0mm, clip]{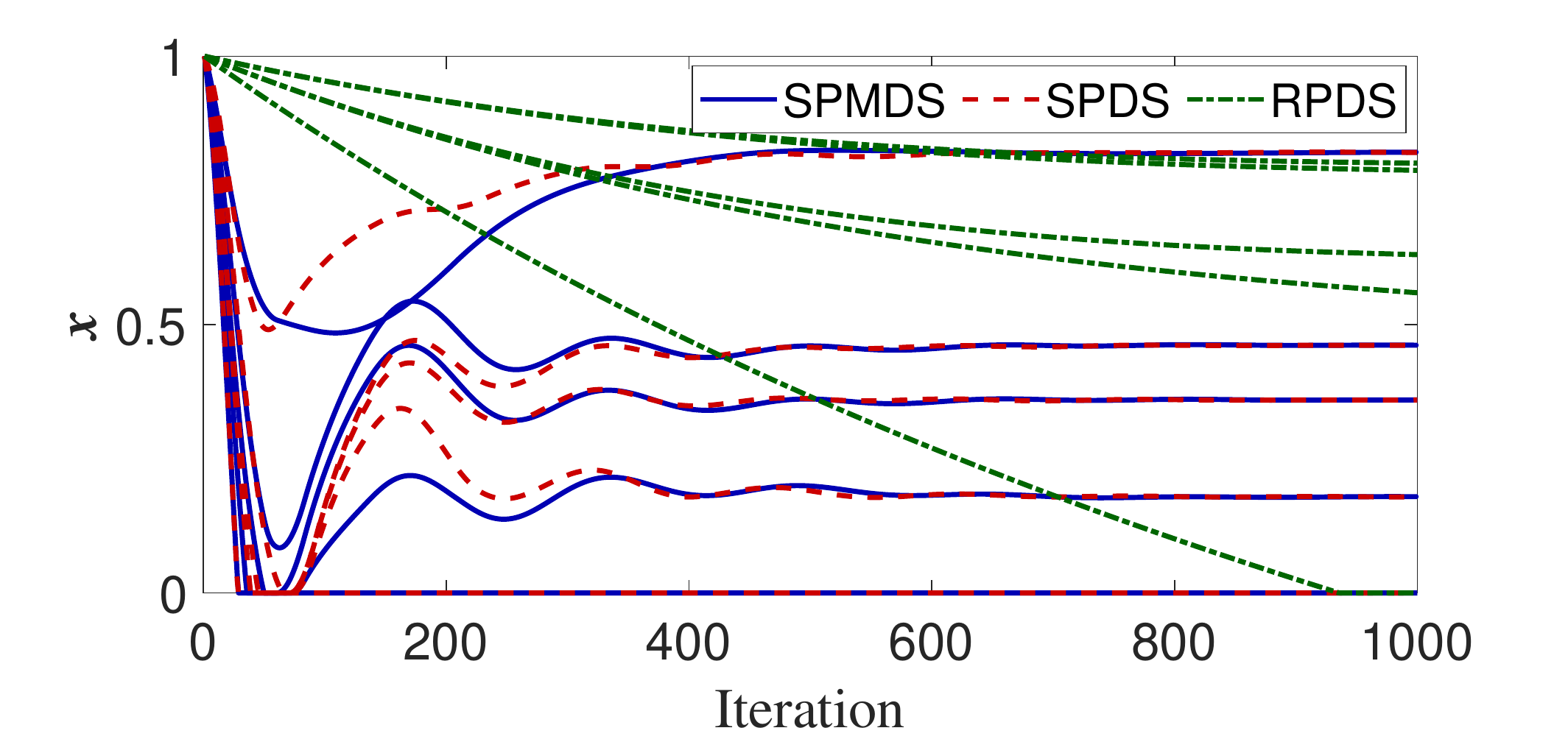}
    \caption{Convergence of $\bm{x}$ for SPMDS, SPDS and RPDS.}
    \label{traffic_x_spmds}
\end{figure}

To compare and analyze the proposed SPMDS with existing decentralized algorithms, we first solve \eqref{32ssss} via a centralized solver to obtain the global optimum $\bm{x}^*$, then we compare the primal convergence errors between RPDS \cite{koshal2011multiuser}, SPDS \cite{liu2017decentralized} and SPMDS for fairness and clarity. The convergence error $\chi^{(\ell)}$ at the $\ell$th iteration is defined as the gap between the solution $\bm{x}^{(\ell)}$ in the $\ell$th iteration and the optimal solution $\bm{x}^*$ as
\begin{equation}
    \chi^{(\ell)} = \|\bm{x}^{(\ell)} - \bm{x}^*\|_2.
\end{equation}
Fig. \ref{traffic_compare} presents the convergence errors of the RPDS, SPDS and SPMDS within 1,000 iterations. The initial points are set as $\bm{1}_5$ uniformly, Fig. \ref{traffic_x_spmds} illustrates the convergence evolution of $\bm{x}$ in RPDS, SPDS and SPMDS, respectively.
For RPDS, all agents use the same regularization parameter $0.1$, and for both the SPDS and SPMDS the primal and dual step sizes are uniformly set as $\alpha = 10^{-3}$ and $\beta = 0.5$. It can be observed in Fig. \ref{traffic_compare} and Fig. \ref{traffic_x_spmds} that both SPDS and SPMDS converge in 600 iterations and the regularization errors are eliminated compared with the RPDS. Moreover, compared with SPDS, SPMDS has the advantage of lower computational cost owing to the dimension reduction strategy while exhibits a similar high convergence speed. Note that the iteration number required to converge may differ in different applications due to network structure and optimization problem formulation.

\section{Conclusion}
In this paper, we have focused on a class of networked optimization problems that minimize convex strongly coupled objective functions under local constraints and globally coupled inequality constraints. A novel decentralized optimization framework with a reduced-dimension primal and multi-dual architecture was proposed to solve this type of problems. In this framework, a dimension reduction technique was developed to virtually group primal variables via $\mathcal{K}$-means and construct subsets of the global constraints correspondingly. The proposed decentralized optimization algorithm, SPMDS, which integrates the dimension reduction technique, is two-facet scalable \emph{w.r.t.}  the agent population size owing to the decentralized architecture, and  \emph{w.r.t.}  the network dimension owing to the primal-multi-dual architecture. The efficiency, efficacy, and convergence of the  proposed decentralized optimization framework was demonstrated through decentralized EV charging control and traffic congestion control problems. Rigorous theoretical analyses of the convergence and optimality have been provided.
\appendices



\section{Proof of Theorem 1}
From \eqref{4sn}, we have the voltage sub-vectors $\bm{\hat{V}}_s(T) \in \mathbb{R}^{n-d} \ \text{for} \ s=1,\ldots,r$. Then in order to satisfy both \eqref{4sn} and \eqref{5sn}, the following inequality should hold
\begin{equation}
    r(n-d) \geq n.
    \label{36sn}
\end{equation}
Therefore, $d \leq n - \frac{n}{r}$, and because of $d \in \mathbb{Z}_{0+}$, the maximum of $d$ should be $d = n - \lceil \frac{n}{r} \rceil$, i.e., as long as \eqref{6sn} stands. This completes the proof.

\section{Proof of Theorem 2}
Define a mapping 
\begin{align}
\Phi(\zeta)\triangleq\left[\begin{array}{c} \tilde{\nabla}_{\mathcal{U}} \mathcal{L}(\mathcal{U}, \bm{\lambda}_e) + \frac{\alpha-\tau_{\mathcal{U}}}{\alpha}\mathcal{U} \\
-\tilde{\nabla}_{\lambda_1} \mathcal{L}(\mathcal{U}, \bm{\lambda}_1)+\frac{\beta-\tau_{\lambda}}{\beta} \bm{\lambda}_1\\
\vdots\\
-\tilde{\nabla}_{\lambda_r} \mathcal{L}(\mathcal{U}, \bm{\lambda}_r)+\frac{\beta-\tau_{\lambda}}{\beta} \bm{\lambda}_r\\
\end{array}\right]
{=}\left[\begin{array}{l}\Phi_{1}(\zeta) \\
\Phi_{21}(\zeta)\\
\vdots\\
\Phi_{2r}(\zeta)\\
\end{array}\right]. \nonumber
\end{align}
Let $\zeta^{*} = \col( \mathcal{U}^{*},\bm{\lambda}_1^{*},\ldots,\bm{\lambda}_r^{*})$ denote the optimizer. By using the decomposable structure of $\mathbb{U}$ and the nonexpansive property of 
$\Pi_{\mathbb{O}}(\bm{o})$, we have
\begin{align}
&\left\|\mathcal{U}^{(\ell+1)}-\mathcal{U}^{*}\right\|_{2}^{2} \nonumber \\ \nonumber
\begin{split} & \quad \leq 
\left\Vert \frac{1}{\tau_{\mathcal{U}}} \Pi_{\mathbb{U}}(\tau_{\mathcal{U}} \mathcal{U}^{(\ell)}-\alpha \tilde{\nabla}_{\mathcal{U}} \mathcal{L}(\mathcal{U}^{(\ell)}, \bm{\lambda}_e^{(\ell)})) \right. \\
& \qquad \left.-\frac{1}{\tau_{\mathcal{U}}} \Pi_{\mathbb{U}}(\tau_{\mathcal{U}} \mathcal{U}^{*}-\alpha \tilde{\nabla}_{\mathcal{U}} \mathcal{L}(\mathcal{U}^{*}, \bm{\lambda}_e^{*})) \right\rVert_{2}^{2} \end{split}\\ \nonumber
\begin{split} & \quad\leq \frac{1}{\tau_{\mathcal{U}}^{2}} \left\Vert \tau_{\mathcal{U}} \mathcal{U}^{(\ell)}-\alpha \tilde{\nabla}_{\mathcal{U}} \mathcal{L}(\mathcal{U}^{(\ell)}, \bm{\lambda}_e^{(\mathcal{\ell})})-\tau_{\mathcal{U}} \mathcal{U}^{*} \right. \\ \nonumber
& \qquad \left. \vphantom{\nabla_{\mathcal{U}} \mathcal{L}\left(\mathcal{U}^{(\ell)}, \lambda^{(\mathcal{\ell})}\right)} +\alpha \tilde{\nabla}_{\mathcal{U}} \mathcal{L}(\mathcal{U}^{*}, \bm{\lambda}_e^{*}) \right\rVert_{2}^{2} \end{split}\nonumber\\
& \quad=\frac{\alpha^2}{\tau_{\mathcal{U}}^{2}}\left\|\mathcal{U}^{(\ell)}-\mathcal{U}^{*}\right\|_{2}^{2}+\frac{\alpha^{2}}{\tau_{\mathcal{U}}^{2}}\|\Phi_{1}(\zeta^{(\ell)})-\Phi_{1}\left(\zeta^{*})\right\|_{2}^{2} \nonumber\\
& \qquad-2 \frac{\alpha^2}{\tau_{\mathcal{U}}^{2}}(\Phi_{1}(\zeta^{(\ell)})-\Phi_{1}(\zeta^{*}))^{\mathsf{T}}(\mathcal{U}^{(\ell)}-\mathcal{U}^{*}).
\label{18sss}
\end{align}

Similarly, for the dual variable $\bm{\lambda}_s$ of group $s$
\begin{align}
&\left\|\bm{\lambda}_s^{(\ell+1)}-\bm{\lambda}_s^{*}\right\|_{2}^{2} \nonumber \nonumber \\ 
\begin{split} & \quad \leq 
\left\Vert \frac{1}{\tau_{\lambda}} \Pi_{\mathbb{D}}\left(\tau_{\lambda} \bm{\lambda}_s^{(\ell)}+\beta \tilde{\nabla}_{\bm{\lambda}_s} \mathcal{L}\left(\mathcal{U}^{(\ell)}, \bm{\lambda}_s^{(\ell)}\right)\right) \right. \\
& \qquad \left.-\frac{1}{\tau_{\lambda}} \Pi_{\mathbb{D}}\left(\tau_{\lambda
} \bm{\lambda}_s^{*}+\beta \tilde{\nabla}_{\bm{\lambda}_s} \mathcal{L}\left(\mathcal{U}^{*}, \bm{\lambda}_s^{*}\right)\right) \right\rVert_{2}^{2} \end{split}\nonumber\\
& \quad\leq \frac{1}{\tau_{\lambda}^{2}} \Big\Vert \tau_{\lambda}(\bm{\lambda}_s^{(\ell)}-\bm{\lambda}_s^{*})+\beta (\tilde{\nabla}_{\bm{\lambda}_s} \mathcal{L}(\mathcal{U}^{(\ell)}, \bm{\lambda}_s^{(\ell)}) \nonumber \\
&\left. \qquad  - \tilde{\nabla}_{\bm{\lambda}_s} \mathcal{L}\left(\mathcal{U}^{*}, \bm{\lambda}_s^{*}\right)) \right\Vert_{2}^{2} \nonumber\\
& \quad=\frac{\beta^2}{\tau_{\lambda}^{2}}\left\|\bm{\lambda}_s^{(\ell)}-\bm{\lambda}_s^{*}\right\|_{2}^{2}+\frac{\beta^{2}}{\tau_{\lambda}^{2}}\|\Phi_{2s}(\zeta^{(\ell)})-\Phi_{2s}\left(\zeta^{*})\right\|_{2}^{2} \nonumber\\
& \qquad-2 \frac{\beta^2}{\tau_{\lambda}^{2}}(\Phi_{2s}(\zeta^{(\ell)})-\Phi_{2s}(\zeta^{*}))^{\mathsf{T}}(\bm{\lambda}_s^{(\ell)}-\bm{\lambda}_s^{*}).
\label{17ss}
\end{align}

\noindent We then can readily have
\begin{align}
&\sum_{s=1}^{r}\left\|\bm{\lambda}_s^{(\ell+1)}-\bm{\lambda}_s^{*}\right\|_{2}^{2} \nonumber \\ 
& \quad \leq \sum_{s=1}^{r}\left(\frac{\beta^2}{\tau_{\lambda}^{2}}\left\Vert \bm{\lambda}_s^{(\ell)}-\bm{\lambda}_s^{*} \right\rVert_{2}^{2}+\frac{\beta^2}{\tau_{\lambda}^{2}}\left\Vert \Phi_{2s}(\zeta^{(\ell)}) - \Phi_{2s}(\zeta^*) \right\rVert_{2}^{2} \right. \nonumber\\
&\left. \qquad ~~~~~~- 2\frac{\beta^2}{\tau_{\lambda}^{2}}(\Phi_{2s}(\zeta^{(\ell)}) - \Phi_{2s}(\zeta^*))^{\mathsf{T}}(\bm{\lambda}_s^{(\ell)}-\bm{\lambda}_s^{*})\right).
\label{11ss}
\end{align}

Let $\mathcal{V}^{(\ell)}$ denote $ \mathcal{V} (\mathcal{U}^{(\ell)},\bm{\lambda}_1^{(\ell)},\ldots,\bm{\lambda}_r^{(\ell)})$, then substitute \eqref{18sss} and \eqref{11ss} into \eqref{6s} gives 

\begin{align}
&\mathcal{V}^{(\ell+1)} \nonumber \\
&\leq \beta^2 \left( \frac{\alpha^2}{\tau_{\mathcal{U}}^{2}}\left\|\mathcal{U}^{(\ell)}-\mathcal{U}^{*}\right\|_{2}^{2}+\frac{\alpha^{2}}{\tau_{\mathcal{U}}^{2}}\|\Phi_{1}(\zeta^{(\ell)})-\Phi_{1}\left(\zeta^{*})\right\|_{2}^{2} \right.\nonumber\\ 
&\quad\left.~~~~~~ -2 \frac{\alpha^2}{\tau_{\mathcal{U}}^{2}}(\Phi_{1}(\zeta^{(\ell)})-\Phi_{1}(\zeta^{*}))^{\mathsf{T}}(\mathcal{U}^{(\ell)}-\mathcal{U}^{*}) \right) \nonumber\\
&\quad+ \alpha^2 \sum_{s=1}^{r} \left(\frac{\beta^2}{\tau_{\lambda}^{2}}\left\Vert \bm{\lambda}_s^{(\ell)}-\bm{\lambda}_s^{*} \right\rVert_{2}^{2} \right.\nonumber\\
&\left.~~~~~~~~~~~~~~~+\frac{\beta^2}{\tau_{\lambda}^{2}}\left\Vert \Phi_{2s}(\zeta^{(\ell)}) - \Phi_{2s}(\zeta^*) \right\rVert_{2}^{2} \right.\nonumber\\
&\left.~~~~~~~~~~~~~~~ - 2\frac{\beta^2}{\tau_{\lambda}^{2}}(\Phi_{2s}(\zeta^{(\ell)}) - \Phi_{2s}(\zeta^*))^{\mathsf{T}}(\bm{\lambda}_s^{(\ell)}-\bm{\lambda}_s^{*})\right). \nonumber\\
&=  \frac{\alpha^2\beta^2}{\tau_{\mathcal{U}}^{2}}\left\|\mathcal{U}^{(\ell)}-\mathcal{U}^{*}\right\|_{2}^{2} +  \frac{\alpha^2\beta^2}{\tau_{\lambda}^{2}} \sum_{s=1}^{r}\left\Vert \bm{\lambda}_s^{(\ell)}-\bm{\lambda}_s^{*} \right\rVert_{2}^{2} \nonumber\\
&\quad+ \frac{\alpha^2\beta^2}{\tau_{\mathcal{U}}^{2}}\|\Phi_{1}(\zeta^{(\ell)})-\Phi_{1}\left(\zeta^{*})\right\|_{2}^{2} \nonumber\\
&\quad+ \frac{\alpha^2\beta^2}{\tau_{\lambda}^{2}}\sum_{s=1}^{r}\left\Vert \Phi_{2s}(\zeta^{(\ell)}) - \Phi_{2s}(\zeta^*) \right\rVert_{2}^{2} \nonumber\\
&\quad-2 \frac{\alpha^2\beta^2}{\tau_{\mathcal{U}}^{2}}(\Phi_{1}(\zeta^{(\ell)})-\Phi_{1}(\zeta^{*}))^{\mathsf{T}}(\mathcal{U}^{(\ell)}-\mathcal{U}^{*})\nonumber\\ \nonumber
\end{align}
\begin{align}
&\quad- 2\frac{\alpha^2\beta^2}{\tau_{\lambda}^{2}}\sum_{s=1}^{r}(\Phi_{2s}(\zeta^{(\ell)}) - \Phi_{2s}(\zeta^*))^{\mathsf{T}}(\bm{\lambda}_s^{(\ell)}-\bm{\lambda}_s^{*}) \nonumber\\
& \leq \mathcal{V}^{(\ell)} +  \left(\frac{\alpha^2\beta^2}{\tau_{\mathcal{U}}^{2}}-\beta^2\right)\left\|\mathcal{U}^{(\ell)}-\mathcal{U}^{*}\right\|_{2}^{2} \nonumber\\
&\quad + \left(\frac{\alpha^2\beta^2}{\tau_{\lambda}^{2}}-\alpha^2\right) \sum_{s=1}^{r}\left\Vert \bm{\lambda}_s^{(\ell)}-\bm{\lambda}_s^{*} \right\rVert_{2}^{2} \nonumber\\
&\quad+ \max \left\{\frac{\alpha^2\beta^2}{\tau_{\mathcal{U}}^{2}},\frac{\alpha^2\beta^2}{\tau_{\lambda}^{2}}\right\}\|\Phi(\zeta^{(\ell)})-\Phi\left(\zeta^{*})\right\|_{2}^{2} \nonumber\\
&\quad-\min \left\{2\frac{\alpha^2\beta^2}{\tau_{\mathcal{U}}^{2}},2\frac{\alpha^2\beta^2}{\tau_{\lambda}^{2}}\right\}(\Phi(\zeta^{(\ell)})-\Phi(\zeta^{*}))^{\mathsf{T}}(\zeta^{(\ell)}-\zeta^{*}).
\label{20ss}
\end{align}

We first deal with the last term of \eqref{20ss} on the right-hand side. Let $d_s(\mathcal{U})$, defined in \eqref{eq:20b}, be partitioned by
\begin{equation}
    d_s(\mathcal{U}) = [d_{s,1}(\mathcal{U}) \ d_{s,2}(\mathcal{U}) \ \cdots \ d_{s,HK}(\mathcal{U}) ]^{\mathsf{T}} \in \mathbb{R}^{HK},
\end{equation}
where $H=n-d$ is the dimension of the voltage subsets,
\begin{equation*}
    d_{s,j}(\mathcal{U}) = \mathcal{Y}_{bs,j} - D_{d,j}^{\mathsf{T}}\mathcal{U}, \quad j=1,\ldots,H,
\end{equation*}
and $D_{d,j}$ is the $j$th row of the matrix $D_{d}$. It can be readily obtained that 
\begin{align}
&(\Phi(\zeta^{(\ell)}) - \Phi(\zeta^*))^{\mathsf{T}}(\zeta^{(\ell)}-\zeta^*)\nonumber\\
& \quad = (\tilde{\nabla}_{\mathcal{U}} \mathcal{L}(\mathcal{U}^{(\ell)}, \bm{\lambda}_e^{(\ell)}) - \tilde{\nabla}_{\mathcal{U}} \mathcal{L}(\mathcal{U}^*, \bm{\lambda}_e^*))^{\mathsf{T}}(\mathcal{U}^{(\ell)}-\mathcal{U}^{*}) \nonumber\\
& \quad \quad + \frac{\alpha-\tau_{\mathcal{U}}}{\alpha}\left\| \mathcal{U}^{(\ell)}-\mathcal{U}^*\right\|_2^2\nonumber\\
&\qquad +\sum_{s=1}^r(-\tilde{\nabla}_{\bm{\lambda}_s} \mathcal{L}(\mathcal{U}^{(\ell)}, \bm{\lambda}_s^{(\ell)})+\tilde{\nabla}_{\bm{\lambda}_s} \mathcal{L}(\mathcal{U}^{*}, \bm{\lambda}_s^{*})
)^{\mathsf{T}}(\bm{\lambda}_s^{(\ell)}-\bm{\lambda}_s^{*}) \nonumber \\
&\quad \quad+ \sum_{s=1}^r \frac{\beta-\tau_{\lambda}}{\beta}\left\|\bm{\lambda}_s^{(\ell)}-\bm{\lambda}_s^{*}\right\|_{2}^{2}.
\label{17s}
\end{align}
Substituting \eqref{modified_gradient} into \eqref{17s}, we have 
\begin{align}
&(\Phi(\zeta^{(\ell)}) - \Phi(\zeta^*))^{\mathsf{T}}(\zeta^{(\ell)}-\zeta^*)\nonumber\\
& \quad = (\tilde{\nabla} G(\mathcal{U}^{(\ell)})- \tilde{\nabla} G(\mathcal{U}^*))^{\mathsf{T}}(\mathcal{U}^{(\ell)} - \mathcal{U}^*)  + \rho\left\| \mathcal{U}^{(\ell)}-\mathcal{U}^*\right\|_2^2\nonumber\\
&\qquad + \frac{\alpha-\tau_{\mathcal{U}}}{\alpha }\left\|\mathcal{U}^{(\ell)}-\mathcal{U}^{*}\right\|_{2}^{2}+\frac{\beta-\tau_{\lambda}}{\beta} \sum_{s=1}^r \left\|\bm{\lambda}_s^{(\ell)}-\bm{\lambda}_s^{*}\right\|_{2}^{2}\nonumber\\
& \qquad + (\tilde{\nabla}_\mathcal{U}d^{\mathsf{T}}(\mathcal{U}^{(\ell)})\bm{\lambda}_e^{(\ell)} - \tilde{\nabla}_\mathcal{U}d^{\mathsf{T}}(\mathcal{U}^*)\bm{\lambda}_e^*))^{\mathsf{T}}(\mathcal{U}^{(\ell)} - \mathcal{U}^*)\nonumber\\
&\qquad + \sum_{s=1}^r(-d_s(\mathcal{U}^{(\ell)})+d_s(\mathcal{U}^{*})
)^{\mathsf{T}}(\bm{\lambda}_s^{(\ell)}-\bm{\lambda}_s^{*}).
\label{22ss}
\end{align}

The convexity of $\nabla G(\mathcal{U})$ indicates
\begin{equation}
    (\nabla G(\mathcal{U}^{(\ell)}) -\nabla G(\mathcal{U}^*))^{\mathsf{T}}(\mathcal{U}^{(\ell)} - \mathcal{U}^*)\geq 0.
    \label{23s}
\end{equation}
The last two terms in \eqref{22ss} can be rewritten as 
\begin{align}
&(\nabla_\mathcal{U}d^{\mathsf{T}}(\mathcal{U}^{(\ell)})\bm{\lambda}_e^{(\ell)} - \nabla_\mathcal{U}d^{\mathsf{T}}(\mathcal{U}^*)\bm{\lambda}_e^*))^{\mathsf{T}}(\mathcal{U}^{(\ell)} - \mathcal{U}^*)\nonumber\\
&+ \sum_{s=1}^r(-d_s(\mathcal{U}^{(\ell)})+d_s(\mathcal{U}^{*})
)^{\mathsf{T}}(\bm{\lambda}_s^{(\ell)}-\bm{\lambda}_s^{*})\nonumber\\
&  = \sum_{s=1}^r\sum_{j=1}^{HK}\bm{\lambda}_{s,j}^{(\ell)}\left(\nabla_\mathcal{U}d_{s,j}^{\mathsf{T}}(\mathcal{U}^{(\ell)})(\mathcal{U}^{(\ell)}-\mathcal{U}^*) - d_{s,j}^{\mathsf{T}}(\mathcal{U}^{(\ell)})\right. \nonumber\\
& \qquad \qquad \qquad\quad+ d_{s,j}^{\mathsf{T}}(\mathcal{U}^{*})\Big)\nonumber\\
& \quad + \sum_{s=1}^r\sum_{j=1}^{HK}\bm{\lambda}_{s,j}^{*}\Big(-\nabla_\mathcal{U}d_{s,j}^{\mathsf{T}}(\mathcal{U}^{*})(\mathcal{U}^{(\ell)}-\mathcal{U}^*) + d_{s,j}^{\mathsf{T}}(\mathcal{U}^{(\ell)}) \nonumber\\
&  \qquad \qquad \qquad \quad - d_{s,j}^{\mathsf{T}}(\mathcal{U}^{*})\Big).
\label{20s}
\end{align}

Eqn. \eqref{20s} comes from \eqref{eq:23b} and the fact that $\nabla_\mathcal{U}d_j(\mathcal{U}_s) = \nabla_\mathcal{U}d_{s,j}(\mathcal{U}_s) = D_{d,j}^{\mathsf{T}}, \ s=1,\cdots,r$.
Since $d_{s,j}(\mathcal{U}), s= 1,\ldots,r, j= 1,\ldots,H$ is convex, it holds that
\begin{equation}
 d_{s,j}(\mathcal{U})-d_{s,j}(\mathcal{V}) - \nabla d_{s,j}(\mathcal{V})^{\mathsf{T}}(\mathcal{U}-\mathcal{V}) \geq 0.
\end{equation}

Hence, the last two terms on the right-hand side of \eqref{22ss} are both non-negative, implying that
\begin{equation}\label{26ss}
\begin{aligned}
&(\Phi(\zeta^{(\ell)}) - \Phi(\zeta^*))^{\mathsf{T}}(\zeta^{(\ell)}-\zeta^*)\\
&\quad \geq F_{\mathcal{U}}\left\| \mathcal{U}^{(\ell)}-\mathcal{U}^*\right\|_2^2 + F_{\lambda} \sum_{s=1}^r \left\|\bm{\lambda}_s^{(\ell)}-\bm{\lambda}_s^{*}\right\|_{2}^{2},
\end{aligned}
\end{equation}

where 
\begin{align}
    F_{\mathcal{U}} &= \rho+\frac{\alpha-\tau_\mathcal{U}}{\alpha}, \nonumber\\
    F_{\lambda} &= \frac{\beta-\tau_\lambda}{\beta}.
\end{align}

To derive an upper bound for the forth term in \eqref{20ss}
\begin{align} \label{23ss}
& \left\Vert \Phi(\zeta^{(\ell)}) - \Phi(\zeta^*) \right\rVert_{2}\nonumber\\
& \quad = \left\Vert \left[
\begin{array}{l}\Phi_{1}(\zeta^{(\ell)})- \Phi_{1}(\zeta^*)\\
\Phi_{21}(\zeta^{(\ell)})-\Phi_{21}(\zeta^*)\\
\vdots\\
\Phi_{2r}(\zeta^{(\ell)}) - \Phi_{2r}(\zeta^*)\\
\end{array}
\right] \right\rVert_2\nonumber\\
& \quad \leq \left\Vert \Phi_{1}(\zeta^{(\ell)})- \Phi_{1}(\zeta^*) \right\rVert_2 + \sum_{s=1}^r \left\Vert \Phi_{2s}(\zeta^{(\ell)})- \Phi_{2s}(\zeta^*) \right\rVert_2 \nonumber\\
 & \quad = \Big\Vert\tilde{\nabla}_{\mathcal{U}} \mathcal{L}(\mathcal{U}^{(\ell)}, \bm{\lambda}_e^{(\ell)}) - \tilde{\nabla}_{\mathcal{U}} \mathcal{L}(\mathcal{U}^*, \bm{\lambda}_e^*)\nonumber\\
  & \quad \qquad +\frac{\alpha-\tau_\mathcal{U}}{\alpha}(\mathcal{U}^{(\ell)}-\mathcal{U}^*)\Big\rVert_2 \nonumber\\
  & \qquad + \sum_{s=1}^r\Big\Vert-\tilde{\nabla}_{\bm{\lambda}_s} \mathcal{L}(\mathcal{U}^{(\ell)}, \bm{\lambda}_s^{(\ell)}) + \tilde{\nabla}_{\bm{\lambda}_s} \mathcal{L}(\mathcal{U}^*, \bm{\lambda}_s^*)\nonumber\\
  & \quad \qquad~~~~~~ +\frac{\beta-\tau_\lambda}{\beta}(\bm{\lambda}_s^{(\ell)}-\bm{\lambda}_s^*)\Big\rVert_2 \nonumber\\ \nonumber
\end{align}
  \begin{align}
  & \quad \leq \left\Vert \tilde{\nabla}_\mathcal{U} G(\mathcal{U}^{(\ell)})-\tilde{\nabla}_\mathcal{U} G(\mathcal{U}^{*}) \right\rVert_2 \nonumber\\
  & \qquad+ (\rho+|\frac{\alpha-\tau_\mathcal{U}}{\alpha}|)\left\Vert \mathcal{U}^{(\ell)}-\mathcal{U}^{*} \right\rVert_2 \nonumber\\
  & \qquad + \left\Vert \tilde{\nabla}_\mathcal{U}d^{\mathsf{T}}(\mathcal{U}^{(\ell)})\bm{\lambda}_e^{(\ell)}-\tilde{\nabla}_\mathcal{U}d^{\mathsf{T}}(\mathcal{U}^{*})\bm{\lambda}_e^{*}\right\rVert_2 \nonumber\\
  & \qquad + \sum_{s=1}^r \left\Vert d_s(\mathcal{U}^{(\ell)})-d_s(\mathcal{U}^{*}) \right\rVert_2  \nonumber \\
  &\qquad+|\frac{\beta-\tau_\lambda}{\beta}|\sum_{s=1}^r \left\Vert \bm{\lambda}_s^{(\ell)}-\bm{\lambda}_s^{*}\right\rVert_2.
\end{align}

The first term on the right-hand side of \eqref{23s} gives 
\begin{equation}
\begin{aligned}\left\|\nabla_{\mathcal{U}} G(\mathcal{U}^{(\ell)})-\nabla_{\mathcal{U}} G(\mathcal{U}^{*})\right\|_{2} &=\left\|\tilde{P}^{\mathsf{T}} \tilde{P}\left(\mathcal{U}^{(\ell)}-\mathcal{U}^{*}\right)\right\|_{2} \\ & \leq L_{\nabla G}\left\|\mathcal{U}^{(\ell)}-\mathcal{U}^{*}\right\|_{2}, \end{aligned}
\end{equation}
where $L_{\nabla G} = nK \max_{i=1,\ldots,H}\{\bar{P}_i^2\}$. Arbitrarily choosing $\mathcal{U}_1$ and $\mathcal{U}_2$ to replace $\mathcal{U}^{(\ell)}$ and $\mathcal{U}^*$, we end up with the Lipschitz continuity of $\nabla_{\mathcal{U}} G(\mathcal{U})$ with the Lipschitz constant $L_{\nabla G}$.

The next step follows the mean-value theorem of
vector-valued functions as follows. 

\textit{\textbf{Theorem 3:} \textbf{(Mean-Value Theorem \cite{apostol1964mathematical})} Let $S \subseteq \mathbb{R}^n$ and the mapping $f:S \mapsto \mathbb{R}^m$ is differentiable at at each point of $S$. Let $\bm{x}$ and $\bm{y}$ be two points in $S$ such that all points between $\bm{x}$ and $\bm{y}$ are in $S$. Then for every vector $\bm{a} \in \mathbb{R}^m$, there is a point $\bm{z}$ between $\bm{x}$ and $\bm{y}$ such that
\begin{equation}
\bm{a}^{\mathsf{T}}(f(\bm{x})-f(\bm{y}))=\bm{a}^{\mathsf{T}}(\nabla f(\bm{z})(f(\bm{x})-f(\bm{y})))
\end{equation}
Further, if $\bm{a}$ is a unit vector such that $\left \Vert \bm{a} \right\rVert_2 = 1$, then it holds that
\begin{equation}
\begin{aligned}\|f(\bm{x})-f(\bm{y})\|_{2} & \leq\|\nabla f(\bm{z})(f(\bm{x})-f(\bm{y}))\|_{2} \\ & \leq L_{f}\|f(\bm{x})-f(\bm{y})\|_{2},
\end{aligned}
\end{equation}
where $\sum_{j=1}^{m}\left\|\nabla f_{j}(\bm{z})\right\|_{2} \leq L_{f}$. \hfill $\blacksquare$} 

Applying \textit{\textbf{Theorem 3}} to the mapping $d_s(\mathcal{U}):\mathbb{R}^{nK} \mapsto \mathbb{R}^{HK}$ yields
\begin{equation}
\left\|d_s(\mathcal{U}^{(\ell)})-d_s(\mathcal{U}^{*})\right\|_{2} \leq L_{d}\left\|\mathcal{U}^{(\ell)}-\mathcal{U}^{*}\right\|_{2},
\end{equation}
where $\sum_{j=1}^{HK}\left\| \nabla_{\mathcal{U}} d_{s,j}(\mathcal{U})\right\|_{2} \leq L_{d}
$. Since $\nabla_{\mathcal{U}} d_{s,j}(\mathcal{U}) = D_{d,j}^{\mathsf{T}}$, we have 
\begin{equation}
\left\|\nabla_{\mathcal{U}} d_{s,j}(\mathcal{U})\right\|_{2} \leq \max _{j=1, \ldots, H K}\left\|D_{d,j}^{\mathsf{T}}\right\|_{2},
\end{equation}
and 
\begin{equation}
    \sum_{j=1}^{HK}\left\|\nabla_{\mathcal{U}} d_{s,j}(\mathcal{U})\right\|_{2} \leq HK\max _{j=1, \ldots, H K}\left\|D_{d,j}^{\mathsf{T}}\right\|_{2} = L_{d}.
\end{equation}
Further, we can readily have 
\begin{equation}
    \sum_{j=1}^{HK}\left\|\nabla_{\mathcal{U}} d_{s,j}(\mathcal{U})\right\|_{2}^2 \leq  \left(  \sum_{j=1}^{HK} \left\|\nabla_{\mathcal{U}} d_{s,j}(\mathcal{U})\right\|_{2} \right)^2 \leq L_d^2,
\end{equation}
This indicates that the mapping $d_{s}(\mathcal{U}): \mathbb{R}^{nK} \mapsto \mathbb{R}^{HK}$ is Lipschitz with the constant $L_d$.
Thus, we have
\begin{equation}
    \sum_{s=1}^r \left\Vert d_s(\mathcal{U}^{(\ell)})-d_s(\mathcal{U}^{*}) \right\rVert_2 \leq rL_d.  
\end{equation}
Consequently, \eqref{23s} becomes 
\begin{align}
& \left\Vert \Phi(\zeta^{(\ell)}) - \Phi(\zeta^*) \right\rVert_{2}\nonumber\\
& \quad \leq (\rho+|\frac{\alpha-\tau_\mathcal{U}}{\alpha}|+L_{\nabla G}+rL_d)\left\Vert \mathcal{U}^{(\ell)}-\mathcal{U}^{*} \right\rVert_2 \nonumber\\
& \qquad + |\frac{\beta-\tau_\lambda}{\beta}|\sum_{s=1}^r \left\Vert \bm{\lambda}_s^{(\ell)}-\bm{\lambda}_s^{*}\right\rVert_2 \nonumber\\
& \qquad + \sum_{s=1}^r \left\Vert \nabla_\mathcal{U}d_s^{\mathsf{T}}(\mathcal{U}^{(\ell)})\bm{\lambda}_s^{(\ell)}-\nabla_\mathcal{U}d_s^{\mathsf{T}}(\mathcal{U}^{*})\bm{\lambda}_s^{*}\right\rVert_2 \nonumber\\
& \quad = L_{\mathcal{U}}\left\Vert \mathcal{U}^{(\ell)}-\mathcal{U}^{*} \right\rVert_2 + |\frac{\beta-\tau_\lambda}{\beta}|\sum_{s=1}^r \left\Vert \bm{\lambda}_s^{(\ell)}-\bm{\lambda}_s^{*}\right\rVert_2 \nonumber\\
& \qquad + \sum_{s=1}^r \left\Vert \sum_{j=1}^{HK} (\bm{\lambda}_s^{(\ell)}-\bm{\lambda}_s^{*})D_{d,j}^{\mathsf{T}}\right\rVert_2 \nonumber\\
& \quad \leq L_{\mathcal{U}}\left\Vert \mathcal{U}^{(\ell)}-\mathcal{U}^{*} \right\rVert_2 +|\frac{\beta-\tau_\lambda}{\beta}|\sum_{s=1}^r \left\Vert \bm{\lambda}_s^{(\ell)}-\bm{\lambda}_s^{*}\right\rVert_2 \nonumber\\
& \qquad + \sum_{s=1}^r  \sum_{j=1}^{HK}\left\Vert \bm{\lambda}_s^{(\ell)}-\bm{\lambda}_s^{*} \right\rVert_2 \left\Vert D_{d,j}^{\mathsf{T}}\right\rVert_2 \nonumber\\
& \quad \leq L_{\mathcal{U}}\left\Vert \mathcal{U}^{(\ell)}-\mathcal{U}^{*} \right\rVert_2 +|\frac{\beta-\tau_\lambda}{\beta}|\sum_{s=1}^r \left\Vert \bm{\lambda}_s^{(\ell)}-\bm{\lambda}_s^{*}\right\rVert_2 \nonumber\\
& \qquad + L_d\sum_{s=1}^r \left\Vert \bm{\lambda}_s^{(\ell)}-\bm{\lambda}_s^{*} \right\rVert_2 \nonumber\\
&\quad = L_{\mathcal{U}}\left\Vert \mathcal{U}^{(\ell)}-\mathcal{U}^{*} \right\rVert_2
+ L_{\lambda}\sum_{s=1}^r \left\Vert \bm{\lambda}_s^{(\ell)}-\bm{\lambda}_s^{*}\right\rVert_2\nonumber\\
&\quad \leq L_{\phi}\left\Vert\zeta^{(\ell)}-\zeta^{*}\right\rVert_2,
\label{34s}
\end{align}
where 
\begin{align}
    L_{\mathcal{U}} &= \rho+|\frac{\alpha-\tau_\mathcal{U}}{\alpha}|+L_{\nabla G}+rL_d, \nonumber\\
    L_{\lambda} &= |\frac{\beta-\tau_\lambda}{\beta}|+L_d, \nonumber\\
    L_{\phi} &= \left\| \left[L_\mathcal{U},  L_{\lambda} \right]^{\mathsf{T}} \right\|_2.
\end{align}
We can then readily obtain that 
\begin{align}
\left\Vert \Phi(\zeta^{(\ell)}) - \Phi(\zeta^*) \right\rVert_{2}^2 \leq L_{\phi}^2\left\Vert\zeta^{(\ell)}-\zeta^{*}\right\rVert_2^2
\label{41ss}
\end{align}

Substituting \eqref{26ss} and \eqref{41ss} into \eqref{20ss}, we have
\begin{align}
&\mathcal{V}^{(\ell+1)} \nonumber \\
& \leq \mathcal{V}^{(\ell)} +  \left(\frac{\alpha^2\beta^2}{\tau_{\mathcal{U}}^{2}}-\beta^2\right)\left\|\mathcal{U}^{(\ell)}-\mathcal{U}^{*}\right\|_{2}^{2} \nonumber\\
&\quad + \left(\frac{\alpha^2\beta^2}{\tau_{\lambda}^{2}}-\alpha^2\right) \sum_{s=1}^{r}\left\Vert \bm{\lambda}_s^{(\ell)}-\bm{\lambda}_s^{*} \right\rVert_{2}^{2} \nonumber\\
&\quad+ \max \left\{\frac{\alpha^2\beta^2}{\tau_{\mathcal{U}}^{2}},\frac{\alpha^2\beta^2}{\tau_{\lambda}^{2}}\right\}\|\Phi(\zeta^{(\ell)})-\Phi\left(\zeta^{*})\right\|_{2}^{2} \nonumber\\
&\quad-\min \left\{2\frac{\alpha^2\beta^2}{\tau_{\mathcal{U}}^{2}},2\frac{\alpha^2\beta^2}{\tau_{\lambda}^{2}}\right\}(\Phi(\zeta^{(\ell)})-\Phi(\zeta^{*}))^{\mathsf{T}}(\zeta^{(\ell)}-\zeta^{*}) \nonumber\\
& \leq \mathcal{V}^{(\ell)} +  \left(\frac{\alpha^2\beta^2}{\tau_{\mathcal{U}}^{2}}-\beta^2\right)\left\|\mathcal{U}^{(\ell)}-\mathcal{U}^{*}\right\|_{2}^{2} \nonumber\\
&\quad + \left(\frac{\alpha^2\beta^2}{\tau_{\lambda}^{2}}-\alpha^2\right) \sum_{s=1}^{r}\left\Vert \bm{\lambda}_s^{(\ell)}-\bm{\lambda}_s^{*} \right\rVert_{2}^{2} \nonumber\\
&\quad+ \max \left\{\frac{\alpha^2\beta^2}{\tau_{\mathcal{U}}^{2}},\frac{\alpha^2\beta^2}{\tau_{\lambda}^{2}}\right\} L_\phi^2\left\|\mathcal{U}^{(\ell)}-\mathcal{U}^{*}\right\|_{2}^{2} \nonumber\\
&\quad+ \max \left\{\frac{\alpha^2\beta^2}{\tau_{\mathcal{U}}^{2}},\frac{\alpha^2\beta^2}{\tau_{\lambda}^{2}}\right\} L_\phi^2\sum_{s=1}^{r}\left\Vert \bm{\lambda}_s^{(\ell)}-\bm{\lambda}_s^{*} \right\rVert_{2}^{2} \nonumber\\
&\quad-\min \left\{2\frac{\alpha^2\beta^2}{\tau_{\mathcal{U}}^{2}},2\frac{\alpha^2\beta^2}{\tau_{\lambda}^{2}}\right\}F_{\mathcal{U}}\left\| \mathcal{U}^{(\ell)}-\mathcal{U}^*\right\|_2^2 \nonumber\\
&\quad-\min \left\{2\frac{\alpha^2\beta^2}{\tau_{\mathcal{U}}^{2}},2\frac{\alpha^2\beta^2}{\tau_{\lambda}^{2}}\right\}F_{\lambda} \sum_{s=1}^r \left\|\bm{\lambda}_s^{(\ell)}-\bm{\lambda}_s^{*}\right\|_{2}^{2} \nonumber\\
&= \mathcal{V}^{(\ell)} +  \left(\frac{\alpha^2\beta^2}{\tau_{\mathcal{U}}^{2}}-\beta^2
+\max \left\{\frac{\alpha^2\beta^2}{\tau_{\mathcal{U}}^{2}},\frac{\alpha^2\beta^2}{\tau_{\lambda}^{2}}\right\}L_\phi^2\right.\nonumber\\&\quad\left.-\min \left\{2\frac{\alpha^2\beta^2}{\tau_{\mathcal{U}}^{2}},2\frac{\alpha^2\beta^2}{\tau_{\lambda}^{2}}\right\}F_{\mathcal{U}}\right)\left\|\mathcal{U}^{(\ell)}-\mathcal{U}^{*}\right\|_{2}^{2} \nonumber\\
&\quad + \left(\frac{\alpha^2\beta^2}{\tau_{\lambda}^{2}}-\alpha^2
+ \max \left\{\frac{\alpha^2\beta^2}{\tau_{\mathcal{U}}^{2}},\frac{\alpha^2\beta^2}{\tau_{\lambda}^{2}}\right\} L_\phi^2 \right.\nonumber\\
&\quad\left. -\min \left\{2\frac{\alpha^2\beta^2}{\tau_{\mathcal{U}}^{2}},2\frac{\alpha^2\beta^2}{\tau_{\lambda}^{2}}\right\}F_{\lambda}
\right) \sum_{s=1}^{r}\left\Vert \bm{\lambda}_s^{(\ell)}-\bm{\lambda}_s^{*} \right\rVert_{2}^{2}.
\label{38s}
\end{align}

Because $\left\|\mathcal{U}^{(\ell)}-\mathcal{U}^{*}\right\|_{2}^{2} \geq 0$ and $\sum_{s=1}^{r}\left\Vert \bm{\lambda}_s^{(\ell)}-\bm{\lambda}_s^{*} \right\rVert_{2}^{2} \geq 0$, to guarantee that $V^{(\ell)}$ is non-increasing in each iteration we require  

\begin{subnumcases}{\label{42sss}} 
 \begin{aligned}
    &\left(\frac{\alpha^2\beta^2}{\tau_{\mathcal{U}}^{2}}-\beta^2
+\max \left\{\frac{\alpha^2\beta^2}{\tau_{\mathcal{U}}^{2}},\frac{\alpha^2\beta^2}{\tau_{\lambda}^{2}}\right\}L_\phi^2\right.\\&\quad\left. \qquad \, -\min \left\{2\frac{\alpha^2\beta^2}{\tau_{\mathcal{U}}^{2}},2\frac{\alpha^2\beta^2}{\tau_{\lambda}^{2}}\right\}F_{\mathcal{U}}\right) < 0 \end{aligned}\label{42sss:a}\\
\begin{aligned}
&\left(\frac{\alpha^2\beta^2}{\tau_{\lambda}^{2}}-\alpha^2
+ \max \left\{\frac{\alpha^2\beta^2}{\tau_{\mathcal{U}}^{2}},\frac{\alpha^2\beta^2}{\tau_{\lambda}^{2}}\right\} L_\phi^2 \right.\\
&\quad\left. \qquad \,  -\min \left\{2\frac{\alpha^2\beta^2}{\tau_{\mathcal{U}}^{2}},2\frac{\alpha^2\beta^2}{\tau_{\lambda}^{2}}\right\}F_{\lambda}
\right) < 0 
\end{aligned}. \label{42sss:b}
\end{subnumcases}

Let $\Psi = \max \left\{\frac{\alpha^2\beta^2}{\tau_{\mathcal{U}}^{2}},\frac{\alpha^2\beta^2}{\tau_{\lambda}^{2}}\right\}$, $\mu \Psi = \min \left\{\frac{\alpha^2\beta^2}{\tau_{\mathcal{U}}^{2}},\frac{\alpha^2\beta^2}{\tau_{\lambda}^{2}}\right\}$ where $0<\mu<1$. Let $M = \frac{\alpha^2\beta^2}{\tau_{\mathcal{U}}^{2}}-\beta^2$ and $N = \frac{\alpha^2\beta^2}{\tau_{\lambda}^{2}}-\alpha^2$, then \eqref{42sss} can be rewritten as 

\begin{subnumcases}{\label{59sss}} 
 \begin{aligned}
    &M+\Psi L_\phi^2  -2\mu \Psi F_{\mathcal{U}} < 0 \end{aligned}\label{59sss:a}\\
\begin{aligned}
& N+ \Psi L_\phi^2 - 2\mu \Psi F_{\lambda} < 0 
\end{aligned}. \label{59sss:b}
\end{subnumcases}
Solving \eqref{59sss} and combing $0<\mu<1$, we have 
\begin{align}
    \max \left\{ \frac{M+\Psi L_\phi^2}{2 \Psi F_{\mathcal{U}}},\frac{N+ \Psi L_\phi^2}{2 \Psi F_{\lambda}}
    \right\} < \mu < 1.
\end{align}
Let $A<0$ denote \eqref{59sss:a} and $B<0$ denote \eqref{59sss:b}, then we can readily have 
\begin{align}
    \mathcal{V}^{(\ell+1)} \leq \mathcal{V}^{(\ell)} + A\left\|\mathcal{U}^{(\ell)}-\mathcal{U}^{*}\right\|_{2}^{2} +B\sum_{s=1}^{r}\left\Vert \bm{\lambda}_s^{(\ell)}-\bm{\lambda}_s^{*} \right\rVert_{2}^{2},
\end{align}
and when $\mathcal{U}^{(\ell)} \neq \mathcal{U}^{*}$, $\mathcal{V}^{(\ell+1)} < \mathcal{V}^{(\ell)}$ holds. Iterating the inequality, we have 
\begin{equation}
    \begin{cases}
      \mathcal{V}^{(1)} \leq \mathcal{V}^{(0)} + A\left\|\mathcal{U}^{(0)}-\mathcal{U}^{*}\right\|_{2}^{2} +B\sum_{s=1}^{r}\| \bm{\lambda}_s^{(0)}-\bm{\lambda}_s^{*} \|_2^2\\
   \mathcal{V}^{(2)} \leq \mathcal{V}^{(1)} + A\left\|\mathcal{U}^{(1)}-\mathcal{U}^{*}\right\|_{2}^{2}+B\sum_{s=1}^{r}\| \bm{\lambda}_s^{(1)}-\bm{\lambda}_s^{*} \|_2^2\\
   \ \ \vdots\\
    \mathcal{V}^{(\ell+1)} \leq \mathcal{V}^{(\ell)} + A\left\|\mathcal{U}^{\ell}-\mathcal{U}^{*}\right\|_{2}^{2}+B\sum_{s=1}^{r}\| \bm{\lambda}_s^{(\ell)}-\bm{\lambda}_s^{*} \|_2^2\\
    \end{cases} 
    \label{41s}
\end{equation}
Summing all inequalities in \eqref{41s}, we have
\begin{align}
    0 \leq \mathcal{V}^{(\ell+1)} &\leq \mathcal{V}^{(0)} + A\sum_{k=0}^{\ell}\left\|\mathcal{U}^{(k)}-\mathcal{U}^{*}\right\|_{2}^{2} \nonumber\\
    &\quad+ B\sum_{k=0}^{\ell}\sum_{s=1}^{r}\| \bm{\lambda}_s^{(k)}-\bm{\lambda}_s^{*} \|_{2}^{2}.
\end{align}

$\mathcal{V}^{(0)}$ is bounded and positive semidefinite leading to 
$\|\mathcal{U}^{(\ell)}-\mathcal{U}^{*}\|_{2}^{2} \to 0$ as  $\ell \to \infty$ and $\sum_{s=1}^{r}\|\bm{\lambda}_s^{(\ell)}-\bm{\lambda}_s^{*} \|_2^2 \to 0$ as  $\ell \to \infty$, besides, the strong duality holds, which implies  $\sum_{s=1}^{r}\|\bm{\lambda}^{(\ell)}-\bm{\lambda}^{*}\|_{2}^{2} \to 0$ as  $\ell \to \infty$. This completes the proof.

\bibliographystyle{IEEEtran}
\bibliography{SPMDS}

\end{document}